# Reduction of points in the group of components of the Néron model of a jacobian

Dino Lorenzini

## Introduction

Let $K$ be a complete field with a discrete valuation. Let $\mathcal{O}_K$ denote the ring of integers of $K$, with maximal ideal $(t)$. Let $k$ be the residue field of $\mathcal{O}_K$, assumed to be algebraically closed of characteristic $p \geq 0$. We shall call a *curve* in this article a smooth proper geometrically connected variety $X/K$ of dimension 1. Let $A/K$ denote the jacobian of $X/K$. Let $P$ and $Q$ be two $K$-rational points of $X$. The divisor of degree zero $P-Q$ defines a $K$-rational point of $A/K$. In this article, we study the reduction of the point $P-Q$ in the Néron model of $A/K$ in terms of the reductions of the points $P$ and $Q$ in a regular model $\mathcal{X}/\mathcal{O}_K$ of $X/K$.

Let $A/K$ be any abelian variety of dimension $g$. Denote by $\mathcal{A}/\mathcal{O}_K$ its Néron model. Recall that the special fiber $\mathcal{A}_k/k$ of $\mathcal{A}/\mathcal{O}_K$ is an extension of a finite abelian group $\Phi_K := \Phi_K(A)$, called the group of components, by a smooth connected group scheme $\mathcal{A}_k^0$, the connected component of zero in $\mathcal{A}_k$. We denote by $\pi : A(K) \to \mathcal{A}_k(k)$ the canonical reduction map. We will often abuse notation and also denote by $\pi$ the composition $A(K) \to \mathcal{A}_k(k) \to \Phi_K$.

In [Lor3], the author introduced two functorial filtrations of the prime-to-$p$ part $\Phi_K^{(p)}$ of the group $\Phi_K$. These filtrations are key in the complete description of all possible groups $\Phi_K^{(p)}$ [Edi]. Filtrations for the full group $\Phi_K$ were later introduced by Bosch and Xarles in [B-X]. An example of a functorial subgroup of $\Phi_K$ occuring in one of the filtrations is the group $\Psi_{K,L}$ described below, where $L/K$ denotes the minimal extension of $K$ such that $A_L/L$ has semistable reduction (see [Des], 5.15). More generally, let $M/K$ be any separable extension. Let $\Phi_M$ denote the group of components of $A_M/M$. The functoriality property of the Néron models induces a map

$$\gamma_{K,M} : \Phi_K \longrightarrow \Phi_M,$$

whose kernel is denoted by $\Psi_{K,M}$.

Given two points $P$ and $Q$ in $X(K)$, it is natural to wonder whether it is possible to predict when the reduction of $P-Q$ in $\Phi_K$ belongs to one of the functorial subgroups mentioned above. This question is not easy since even deciding whether the reduction of $P-Q$ is trivial is in general not obvious. We give in this paper a sufficient condition on the special fiber of a model $\mathcal{X}$ for the image of the point $P-Q$ in $\Phi_K$ to belong to the subgroup $\Psi_{K,L}$. When this condition is satisfied, we are able to provide a formula for the order of this image. We conjecture that the sufficient condition alluded to above is also necessary and we provide evidence in support of this conjecture. We also discuss cases where the image of the point $P-Q$ belongs to the subgroup $\Theta_K^{[3]}$ of $\Psi_{K,L}$ (notation recalled in 6.6), using a pairing associated to the group of component $\Phi_K$.

# 1 The main results

Let $X/K$ be a curve. Let $\mathcal{X}/\mathcal{O}_K$ be a regular model of $X/K$. Let $\mathcal{X}_k := \sum_{i=1}^{v} r_i C_i$ denote the special fiber of $\mathcal{X}$ and let $M := ((C_i \cdot C_j))_{1 \leq i,j \leq v}$ be the associated *intersection matrix*. The *dual graph* $G$ associated to $\mathcal{X}_k$ is defined as follows. The vertices of $G$ are the curves $C_i$ and, when $j \neq h$, the vertex $C_j$ is linked in $G$ to the vertex $C_h$ by exactly $(C_j \cdot C_h)$ edges. The degree of the vertex $C_i$ in $G$ is the integer $d_i := \sum_{i \neq j} (C_i \cdot C_j)$.

Let ${}^tR := (r_1, \ldots, r_v)$, so that $MR = 0$. We assume in this paper that $\gcd(r_1, \ldots, r_v) = 1$. The triple $(G, M, R)$ is an example of what we called an *arithmetical graph* in [Lor1]. When the coefficients of $M$ are not thought of as intersection numbers, we may denote $(C_i \cdot C_j)$ simply by $c_{ij}$. As we will recall in section 6, Raynaud has shown that the group of components $\Phi_K(\mathrm{Jac}(X))$ is isomorphic to $\mathrm{Ker}({}^tR)/\mathrm{Im}(M)$, where ${}^tR : \mathbb{Z}^v \to \mathbb{Z}$ and $M : \mathbb{Z}^v \to \mathbb{Z}^v$ are the linear maps associated with the matrices $M$ and ${}^tR$. We call the group $\Phi(G) := \mathrm{Ker}({}^tR)/\mathrm{Im}(M)$ the *group of components* of the arithmetical graph $(G, M, R)$.

Let $(C,r)$ and $(C',r')$ be two vertices of $G$. Let $E(C,C')$ denote the vector of $\mathbb{Z}^v$ with null components everywhere except for $r'/\gcd(r,r')$ in the $C$-component, and $-r/\gcd(r,r')$ in the $C'$-component. Clearly, $E(C,C') \in \mathrm{Ker}({}^tR)$. The image of $E(C,C')$ in the quotient $\mathrm{Ker}({}^tR)/\mathrm{Im}(M)$ will be called the element of $\Phi(G)$ associated to the pair of vertices $(C,C')$.

Let $P \in X(K)$. Let $\overline{P} \in \mathcal{X}$ denote the closure of $P$ in $\mathcal{X}$. The Cartier divisor $\overline{P}$ intersects $\mathcal{X}_k$ in a smooth point of $\mathcal{X}_k$. Hence, there exists a unique component $C_P$ of $\mathcal{X}_k$, of multiplicity one, such that $\overline{P} \cap \mathcal{X}_k \in C_P$. To determine the image of $P-Q$ in $\Phi_K(\mathrm{Jac}(X))$, it is sufficient to determine the image of the vector $E(C_P,C_Q)$ in $\mathrm{Ker}({}^tR)/\mathrm{Im}(M)$. It follows immediately from this description that the image of $P-Q$ in $\Phi_K$ is trivial if $C_P = C_Q$. Thus, in the remainder of this article, we shall usually assume that $C_P \neq C_Q$.

**1.1** Let us recall the following terminology. A *node* of a graph $G$ is a vertex of degree greater than 2. A *terminal vertex* is a vertex of degree 1. The topological space obtained from $G$ by removing all its nodes is the union of connected components. A *chain* of $G$ is a connected subgraph of the closure of such a connected component. In particular, a chain contains at most two nodes of $G$. If a chain contains a terminal vertex, we call it a *terminal chain.* We define the *weight* of a chain $\mathcal{C}$ to be the integer $w(\mathcal{C}) := \gcd(r_j, C_j \text{ a vertex on } \mathcal{C})$. Let $(C,r),(C_1,r_1),\ldots,(C_n,r_n),(C',r')$, be the vertices on a chain $\mathcal{C}$ of an arithmetical graph, with $C$ and $C'$ nodes: then $(C\cdot C_1) = (C_i \cdot C_{i+1}) = (C_n \cdot C') = 1$. The reader will check that $\gcd(r,r_1) = \gcd(r_1,r_2) = \ldots = \gcd(r_n,r')$. In particular, $w(\mathcal{C}) = \gcd(r,r_1)$. When $(C,r),(C_1,r_1),\ldots,(C_n,r_n)$ are the vertices of a terminal chain, with $C_n$ the terminal vertex, then $\gcd(r,r_1) = r_n$. Note that if the set of vertices on a chain consists of exactly two nodes $C$ and $C'$, it may happen that $(C\cdot C') > 1$.

Let $(C,r)$ and $(C',r')$ be two distinct vertices of $G$. We say that the pair $(C,C')$ is *weakly connected* if there exists a path $\mathcal{P}$ in $G$ between $C$ and $C'$ such that, for each edge $e$ on $\mathcal{P}$, the graph $G\setminus\{e\}$ is disconnected. Note that when a pair $(C,C')$ is weakly connected, then the path $\mathcal{P}$ is the unique shortest path between $C$ and $C'$. If a pair is not weakly connected, we will say that it is *multiply connected.* A graph is a tree if and only if every pair of vertices of $G$ is weakly connected.

Let $(C,r)$ and $(C',r')$ be a weakly connected pair with associated path $\mathcal{P}$. While walking on $\mathcal{P}\setminus\{C,C'\}$ from $C$ to $C'$, label each encountered node consecutively by $(C_1, r_1), (C_2, r_2), \ldots, (C_s, r_s)$. (There may be no such nodes, in which case the integer $s$ is set to be 0.) Thus $\mathcal{P}$ is the union of chains: the chain $\mathcal{C}_0$ from $C$ to $C_1$, then the chain $\mathcal{C}_1$ from $C_1$ to $C_2$, and so on. The last chain on $\mathcal{P}$ is the chain $\mathcal{C}_s$ from $C_s$ to $C'$. If there are no nodes on $\mathcal{P} \setminus \{C, C'\}$, then $\mathcal{P}$ is a chain from $C$ to $C'$, and if there are no vertices on $\mathcal{P} \setminus \{C, C'\}$, then by definition of weakly connected, $(C \cdot C') = 1$. Let $\ell$ be a prime number. We say that the weakly connected pair $(C, C')$ is $\ell$-*breakable* if, for all $i = 0, \ldots, s$, the weight $w(\mathcal{C}_i)$ is not divisible by $\ell$. In particular, if the pair $(C, C')$ is $\ell$-breakable, then each chain $\mathcal{C}_i$ contains a vertex of multiplicity prime to $\ell$. To study the element of $\Phi(G)$ associated to the pair $(C, C')$, we will break the graph $G$ at each such vertex and study each smaller graph so obtained individually.

Note that there is only one reduction type of curve of genus $g = 1$ which contains a weakly connected pair that is not $\ell$-breakable: the type $I_\nu^*$, with $\ell = 2$ and $\nu > 0$. For examples with $g > 1$, see 7.6.

**1.2** Let $(C, C')$ be a weakly connected and $\ell$-breakable pair. Let $\mathcal{P}$ denote the associated path between $C$ and $C'$, with nodes $(C_1, r_1), \ldots, (C_s, r_s)$. If $s = 0$, set $\lambda(C, C') := 1$. If $s > 0$, define $\lambda(C, C')$ as follows. Remove all edges of $\mathcal{P}$ from $G$ to obtain a disconnected graph $\mathcal{G}$. Let $\mathcal{G}_i$, $i = 1, \ldots, m$, denote the connected components of $\mathcal{G}$. Let us number these connected components in such a way that the node $(C_i, r_i)$ on the path $\mathcal{P}$ belongs to the graph $\mathcal{G}_i$. Let

$$m_i := \gcd(r_j, (D_j, r_j) \text{ a vertex of } \mathcal{G}_i)$$

and let $\lambda(C, C')$ denote the power of $\ell$ such that

$$\operatorname{ord}_\ell(\lambda(C, C')) := \max\{\operatorname{ord}_\ell(r_i/m_i), C_i \text{ a node on } \mathcal{P}\}.$$

**1.3** Recall that a finite abelian group $H$ can be written as a product $H \cong \prod_{\ell \text{ prime}} H_\ell$. The group $H_\ell$ is called the $\ell$-*part of* $H$. Let $h$ be an element of $H$ of order $m$. We call the $\ell$-*part of* $h$ the following element $h_\ell$ of $H$. If $\ell \nmid m$, then $h_\ell$ is trivial. Otherwise, write $1 = \sum_{\ell \text{ prime}} a_\ell m/\ell^{\operatorname{ord}_\ell(m)}$. Then set $h_\ell := h^{a_\ell m/\ell^{\operatorname{ord}_\ell(m)}}$. The reader will check that the element $h_\ell$ does not depend on the choice of the coefficients $a_\ell$. We may now state the main results of this article.

**Theorem 6.5.** *Let $X/K$ be a curve. Let $\mathcal{X}/\mathcal{O}_K$ be a regular model of $X/K$ with associated arithmetical graph $(G, M, R)$. Let $\ell \neq p$ be a prime. Let*

*$P, Q \in X(K)$ with $C_P \neq C_Q$. If the pair $(C_P, C_Q)$ is weakly connected and $\ell$-breakable, then the image of the $\ell$-part of $P-Q$ in $\Phi_K(\mathrm{Jac}(X))$ belongs to $\Psi_{K,L}$, and has order $\lambda(C_P, C_Q)$.*

**Theorem 7.3/7.4.** *Let $X/K$ be a curve. Let $\mathcal{X}/\mathcal{O}_K$ be a regular model of $X/K$ with associated arithmetical graph $(G, M, R)$. Let $P, Q \in X(K)$ with $C_P \neq C_Q$. If the pair $(C_P, C_Q)$ is not weakly connected, or if it is weakly connected but not $\ell$-breakable for some prime $\ell \neq p$, then the image of $P-Q$ in $\Phi_K(\mathrm{Jac}(X))$ does not belong to $\Psi_{K,L}$.*

Note that Theorem 7.3 is only a partial converse to Theorem 6.5 since 7.3 provides information only on the image of $P-Q$ and not on the image of the $\ell$-part of $P-Q$.

**1.4** Recall that the connected component $\mathcal{A}_k^0$ of the Néron model $\mathcal{A}/\mathcal{O}_K$ is the extension of an abelian variety of dimension $a_K$ by the product of a torus and an unipotent group of dimension $t_K$ and $u_K$ respectively. The integers $a_K$, $t_K$, and $u_K$ are called the abelian, toric, and unipotent ranks of $A/K$, respectively. For each prime $\ell$ dividing $[L:K]$, $\ell \neq p$, let $K_\ell/K$ denote the unique subfield of $L$ with the property that $[K_\ell : K] = \ell^{\mathrm{ord}_\ell([L:K])}$. An abelian variety has *potentially good reduction* if $t_L = 0$. It is said to have *potentially good $\ell$-reduction* if $t_{K_\ell} = 0$. An abelian variety with potentially good reduction has potentially good $\ell$-reduction for all primes $\ell \neq p$, but the converse is false, even when $p = 0$.

We shall say that an element $h$ of a group $H$ is *divisible by $\ell$*, or is *$\ell$-divisible* if there exists $g \in H$ such that $\ell g = h$. Note that the $\ell$-part $h_\ell$ of $h$ is $\ell$-divisible if and only if $h$ is $\ell$-divisible.

**Theorem 8.2.** *Let $X/K$ be a curve. Let $\ell \neq p$ be a prime. Let $P, Q \in X(K)$. Assume that $\mathrm{Jac}(X)/K$ has potentially good $\ell$-reduction. Then $\Psi_{K,L,\ell} = \Phi_K(\mathrm{Jac}(X))_\ell$. If the $\ell$-part of the image of $P-Q$ in $\Phi_K(\mathrm{Jac}(X))$ is not trivial, then $P-Q$ is not divisible by $\ell$ in $\mathrm{Jac}(X)(K)$.*

This article will proceed as follows. In the next four sections, we prove several propositions on arithmetical graphs needed to compute the order in $\Phi_K$ of elements of the form $\pi(P-Q)$. In particular, we introduce in the third section a very useful pairing on $\Phi \times \Phi$ that is non-degenerate. These first four sections are linear algebraic in essence and can be read independently of the rest of the paper. In the sixth section, we prove the first theorem stated above. In section seven, we discuss a partial converse to this theorem. In the last section, we study the case where the jacobian has potentially good $\ell$-reduction and prove Theorem 8.2.

# 2 Terminal chains

Let $(G, M, R)$ be an arithmetical graph. As the reader may have noted, it is not easy in general to compute the order of the group $\Phi(G)$, or the order in $\Phi(G)$ of a given pair of vertices of $G$. There is no easy criteria to determine in terms of $G$ whether, for instance, $|\Phi(G)| = 1$ (see, however, 7.5 and 3.3). When the arithmetical graph is reduced, that is, when all its multiplicities are equal to 1, such a criterion exists: $\Phi(G)$ is trivial if and only if $G$ is a tree. We provide in this section a necessary condition for a pair $(C, C')$ to have order 1. When the arithmetical graph is reduced, a necessary and sufficient criterion already exists. Indeed, it is shown in [Lor4], 2.3, that:

**Proposition 2.1** *When $G$ is reduced, a pair has order $1$ if and only if it is weakly connected.*

We shall see below that even in the general case, it is possible to show that certain weakly connected pairs have order 1. After a series of preliminary lemmas on chains, we prove in 2.7 the main result of this section, that $E(C, C')$ is trivial if $C$ and $C'$ both belong to the same terminal chain. The case where $C$ and $C'$ are consecutive vertices on a chain is easy and is treated in the following lemma.

**Lemma 2.2** *Let $(C, r)$ and $(C', r')$ be two vertices of an arithmetical graph $(G, M, R)$ joined by a single edge $e$. Assume that $G \setminus \{e\}$ is disconnected. Let $G_C$ denote the connected component of $G \setminus \{e\}$ that contains $C$. Let $s := \gcd(d, (D, d)$ vertex on $G_C)$. Then the image of $E(C, C')$ in $\Phi(G)$ is killed by $\gcd(r, r')/s$. In particular, if $C$ and $C'$ belong to the same terminal chain, then the image of $E(C, C')$ is trivial.*

*Proof:* Multiply each column of $M$ corresponding to a vertex $(D, d)$ of $G_C$ by $d/s$. Add all these columns to the $C$-column multiplied by $r/s$. The new matrix has the vector $(-\gcd(r, r')/s)E(C, C')$ in the $C$-column. Hence, $(-\gcd(r, r')/s)E(C, C')$ belongs to $\mathrm{Im}(M)$ and is thus trivial in $\Phi(G)$. If $C$ and $C'$ belong to the same terminal chain, we may without loss of generality assume that $G_C$ contains the terminal vertex of the chain. The terminal vertex has then multiplicity $s$, which equals $\gcd(r, r')$.

**2.3** Let $n \geq 1$. Let $(C, r), (C_1, r_1), \ldots, (C_n, r_n), (C', r')$, be the vertices on a chain of an arithmetical graph. Letting $-c_i$ denote the self-intersection of $C_i$, we obtain a $(n \times n)$ matrix $N$ and a relation:

$$N := \begin{pmatrix} -c_1 & 1 & 0 & \ldots & 0 \\ 1 & -c_2 & 1 & & \vdots \\ 0 & 1 & \ddots & \ddots & 0 \\ \vdots & & \ddots & -c_{n-1} & 1 \\ 0 & \ldots & 0 & 1 & -c_n \end{pmatrix} \quad \text{and} \quad N \begin{pmatrix} r_1 \\ r_2 \\ \vdots \\ r_{n-1} \\ r_n \end{pmatrix} = \begin{pmatrix} -r \\ 0 \\ \vdots \\ 0 \\ -r' \end{pmatrix}.$$

It is possible to find a sequence of integers $b_1 = 1, b_2, \ldots, b_n$ such that

$$(b_1, \ldots, b_n) \cdot N = (0, \ldots, 0, -b)$$

for some $b \in \mathbb{Z}$. Indeed, set $b_1 = 1$ and solve for $b_2$ in the above equation. Once $b_1$ and $b_2$ are known, then it is possible to solve for $b_3$, and so on. Let

$$A := \begin{pmatrix} b_1 & b_2 & b_3 & \ldots & b_n \\ 0 & 1 & 0 & \ldots & 0 \\ 0 & 0 & 1 & & \vdots \\ \vdots & & & \ddots & 0 \\ 0 & \ldots & \ldots & 0 & 1 \end{pmatrix} \quad \text{and} \quad \mathcal{R} := \begin{pmatrix} 1 & 0 & \ldots & 0 & r_1 \\ 0 & 1 & & & r_2 \\ \vdots & & \ddots & & \vdots \\ 0 & & & 1 & \vdots \\ 0 & \ldots & \ldots & 0 & r_n \end{pmatrix}.$$

**Lemma 2.4** *We have $br_n = r + b_n r'$. When the chain $(C, r), \ldots, (C_n, r_n)$ is terminal with terminal vertex $C_n$, then $br_n = r$.*

*Proof:* Compute $(AN)\mathcal{R}$ and $A(N\mathcal{R})$ and identify the top right coefficients of these matrices.

Let us note there that the integers $b_1 = 1, b_2, \ldots, b_n, b$, are all positive. Indeed, if $b \leq 0$, then $br_n = r + b_n r' \leq 0$ implies $b_n < 0$. If $b_i < 0$ for some $i$, then the equality $b_i r_{i-1} = r + b_{i-1} r_i$, implies that $b_{i-1} < 0$, which is a contradiction since $b_1 > 0$.

**2.5** The sequence $(C_1, r_1), (C_2, r_2), \dots, (C_n, r_n), (C', r')$, is also a chain, with associated matrix $N^{11}$, the principal minor of $N$ obtained by removing the first row and first column of $N$. Let $d_1 = 1, d_2, \dots, d_{n-1},\ d,$ denote the integers associated to $N^{11}$ such that

$$(d_1, d_2, \dots, d_{n-1}) N^{11} = (0, \dots, 0, -d).$$

Let

$$A' := \begin{pmatrix} -1 & 0 & d_1 & d_2 & \dots & d_{n-1} \\ 0 & b_1 & b_2 & b_3 & \dots & b_n \\ \vdots & 0 & 1 & 0 & \dots & 0 \\ \vdots & 0 & 0 & 1 & & \vdots \\ & \vdots & & & \ddots & 0 \\ 0 & 0 & \dots & \dots & 0 & 1 \end{pmatrix} \text{ and } N' := \begin{pmatrix} \begin{array}{ccccc} 1 & 0 & 0 & \dots & 0 \\ \hline \\ & & \mathrm{N} & & \\ \\ \end{array} \end{pmatrix}.$$

The matrix $A'N'$ is an $(n+1) \times n$ matrix. Using operations involving only the columns of $A'N'$, it is easy to see that $A'N'$ is equivalent over $\mathbb{Z}$ to the following matrix (we shall say that $A'N'$ is 'column equivalent' to):

$$\begin{pmatrix} 0 & 0 & \dots & 0 & -d \\ 0 & 0 & \dots & 0 & -b \\ 1 & 0 & 0 & \dots & 0 \\ 0 & \ddots & \ddots & \ddots & \vdots \\ \vdots & & 1 & 0 & 0 \\ 0 & \dots & 0 & 1 & 0 \end{pmatrix}.$$

Set $d_0 = 0$.

**Lemma 2.6** *Let $(C, r), (C_1, r_1), \dots, (C_n, r_n)$ be a terminal chain of an arithmetical graph. Then $\det(N) = (-1)^n r/r_n$. Moreover, $r$ divides $r_i b_j - b_i r_j$, for all $i \neq j$, and*

$$\frac{r}{r_1} = \frac{r_i b_j - b_i r_j}{r_i d_{j-1} - d_{i-1} r_j}, \text{ for all } \ i \neq j, \quad 1 \leq i, j \leq n.$$

*In particular, $b_n r_1 / r_n$ is congruent to $1$ modulo $r/r_n$ and $\gcd(b_n, r/r_n) = 1$.*

*Proof:* Recall that, with the notation introduced above, we have

$$\begin{aligned} (r_1, \dots, r_n)N &= (-r, 0, \dots, 0), \\ (b_1, \dots, b_n)N &= (0, \dots, 0, -b), \\ (0, d_1, \dots, d_{n-1})N &= (d_1, 0, \dots, 0, -d). \end{aligned}$$

Recall also that $b_1 = d_1 = 1$, and that since the vertices form a terminal chain, Lemma 2.4 shows that $b = r/r_n$ and $d = r_1/r_n$. It is easy to check that $r_n = \gcd(r, r_1)$ and that $r_n$ divides all $r_i$.

Let $N^* := ((a_{ij}))_{1 \le i,j \le n}$ denote the comatrix of $N$: $N^*N = NN^* = \det(N)I_n$. Multiply both sides of the three equalities above by $((a_{ij}))$. We find that

$$\begin{aligned} \det(N)r_i &= -a_{i,1}r \qquad \forall\, i = 1, \dots, n, \\ \det(N)b_i &= -a_{i,n}b \qquad \forall\, i = 1, \dots, n, \\ \det(N)d_i &= a_{i+1,1} - a_{i+1,n}d \qquad \forall\, i = 0, \dots, n-1. \end{aligned}$$

In particular, $\det(N)r_n = -a_{n,1}r = (-1)^n r$. It follows from the three equalities above that

$$\begin{aligned} r_i b_j - r_j b_i &= r_n(a_{i,1}a_{j,n} - a_{j,1}a_{i,n}), \\ (r_i d_{j-1} - r_j d_{i-1})r &= r_n(a_{i,1}a_{j,n} - a_{j,1}a_{i,n})r_1. \end{aligned}$$

From the equality $(r_i b_j - r_j b_i)r_1 = (r_i d_{j-1} - r_j d_{i-1})r$, we find that $r \mid r_i b_j - r_j b_i$. This concludes the proof of Lemma 2.6.

As a corollary to our study of the properties of the matrix $N$, we may now prove the following result.

**Proposition 2.7** *Let $(G, M, R)$ be an arithmetical graph, and let $(C, r)$, $(C_1, r_1), \dots, (C_n, r_n)$ be a terminal chain of $G$. Then $E(C_i, C_j)$ is trivial in $\Phi(G)$, for all $i, j \in \{1, \dots, n\}$, $i \ne j$.*

*Proof:* The matrix $M$ has the form

$$M = \left(\begin{array}{ccc|c} & & * & \\ & & \vdots & \\ * & \cdots & * & 1 \\ \hline & & 1 & \\ & & & N \\ & & & \end{array}\right). \quad \text{Let } A'' := \left(\begin{array}{c|c} \mathrm{Id}_s & \\ \hline & A' \\ & \end{array}\right),$$

where $A'$ is as in 2.5 and, if $v$ denotes the number of vertices of $G$, then $s := v - n - 1$. Then, using 2.5 and the facts that $d = r_1/r_n$ and $b = r/r_n$, the reader will check that $A''M$ is column equivalent to a matrix of the form

$$\left(\begin{array}{cc|cccc} * & * & & & & \\ * & * & 0 & \dots & 0 & -r_1/r_n \\ \hline & 1 & 0 & \dots & 0 & -r/r_n \\ & & 1 & \ddots & & 0 \\ & & & \ddots & 0 & \vdots \\ & & & & 1 & 0 \end{array}\right).$$

The transpose of the vector $A''E(C_i, C_j)$ has the form (if $i < j$):

$$\frac{1}{r_n}(0, \dots, 0, d_{i-1}r_j - d_{j-1}r_i, b_i r_j - r_i b_j, 0, \dots, 0, +r_j, 0, \dots, 0, -r_i, 0, \dots, 0),$$

(where the first $s$ coefficients are 0). We claim that $A''E(C_i, C_j)$ is in the span of the last $n$ columns of the matrix $A''M$. To prove this claim, it is sufficient to show that $\frac{1}{r_n}(d_{i-1}r_j - d_{j-1}r_i, b_i r_j - r_i b_j)$ is an integer multiple of $(-r_1/r_n, -r/r_n)$, which follows immediately from Lemma 2.6. Since $A''$ is invertible over $\mathbb{Z}$, $A''E(C_i, C_j)$ is in the span of $A''M$ if and only if $E(C_i, C_j)$ is in the span of $M$. Hence, $E(C_i, C_j)$ is trivial in $\Phi(G)$.

We conclude this section with a key lemma used in the next sections.

**Lemma 2.8** *Let $(C, r), (C_1, r_1), \dots, (C_n, r_n)$ be a terminal chain. Then*

$$\frac{1}{rr_1} + \frac{1}{r_1 r_2} + \dots + \frac{1}{r_{n-1}r_n} = \frac{b_n}{rr_n}.$$

*Proof:* We proceed by induction on $n$. If $n = 1$, Lemma 2.8 holds since $b_1 = 1$. By induction hypothesis applied to $C_1, \dots, C_n$,

$$\frac{1}{r_1 r_2} + \dots + \frac{1}{r_{n-1}r_n} = \frac{d_{n-1}}{r_1 r_n}.$$

Lemma 2.6 shows that $r/r_1 = (r_1 b_n - r_n b_1)/(r_1 d_{n-1} - d_0 r_n)$. In other words, $d_{n-1}r = r_1 b_n - r_n$. Dividing both sides by $rr_1 r_n$ shows that

$$\frac{d_{n-1}}{r_1 r_n} = \frac{b_n}{rr_n} - \frac{1}{rr_1}.$$

# 3 A pairing attached to $\Phi$

**3.1** Let us introduce in this section a pairing associated to $\Phi(G)$. Let $(G, M, R)$ be any arithmetical graph. Let $\tau, \tau' \in \Phi$ and let $T, T' \in \mathrm{Ker}({}^tR)$ be vectors whose images in $\Phi$ are $\tau$ and $\tau'$, respectively. Let $S, S' \in \mathbb{Z}^v$ be such that $MS = nT$ and $MS' = n'T'$. Note that $n$ and $n'$ are divisible by the order of $\tau$ and $\tau'$, respectively. Define

$$\begin{array}{rcl} \langle\ ;\ \rangle : \Phi \times \Phi & \longrightarrow & \mathbb{Q}/\mathbb{Z} \\ (\tau, \tau') & \mapsto & ({}^tS/n)M(S'/n') \pmod{\mathbb{Z}}. \end{array}$$

It is shown in [Lor5] that this pairing is well-defined and non-degenerate. Moreover, let $(C, r)$ and $(C', r')$ be a weakly connected pair with associated path $\mathcal{P}$. While walking on $\mathcal{P} \setminus \{C, C'\}$ from $C$ to $C'$, label each encountered vertex consecutively by $(C_1, r_1), (C_2, r_2), \ldots, (C_n, r_n)$. The following proposition is proved in [Lor5], XX.

**Proposition 3.2** *Keep the notation introduced above. Assume that $(C, C')$ is a weakly connected pair of $G$. Let $\gamma$ denote the image of the element $E(C, C')$ in $\Phi(G)$. If $(D, s)$ and $(D', s')$ are any two distinct vertices on $G$, let $\delta$ denote the image of $E(D, D')$ in $\Phi(G)$. Let $C_\alpha$ denote the vertex of $\mathcal{P}$ closest to $D$ in $G$, and let $C_\beta$ denote the vertex of $\mathcal{P}$ closest to $D'$. Assume that $\alpha \leq \beta$. (Note that we may have $\alpha = \beta$, and we may have $D = C_\alpha$ or $D' = C_\beta$.) Then*

$$\langle \gamma, \delta \rangle = \mathrm{lcm}(r, r')\mathrm{lcm}(s, s')(1/r_\alpha r_{\alpha+1} + 1/r_{\alpha+1} r_{\alpha+2} + \ldots + 1/r_{\beta-1} r_\beta).$$

*In particular, if $C_\alpha = C_\beta$, then $\langle \gamma, \delta \rangle = 0$. Moreover,*

$$\langle \gamma, \gamma \rangle = \mathrm{lcm}(r, r')^2 (1/r r_1 + 1/r_1 r_2 + \ldots + 1/r_n r').$$

The existence of this perfect pairing has the following interesting consequences.

**Proposition 3.3** *Let $(G, M, R)$ be any arithmetical graph. Let $\ell$ be any prime. Let $(C, C')$ be a weakly connected pair of $(G, M, R)$ such that $\ell \nmid rr'$. If $(C, C')$ is not $\ell$-breakable, then the $\ell$-part of $E(C, C')$ is not trivial in $\Phi_K$.*

*Proof:* By definition of weakly connected and not $\ell$-breakable (1.1), the path $\mathcal{P}$ in $G$ linking $C$ and $C'$ contains two consecutive vertices $(D, d)$ and $(D', d')$ such that $\mathrm{ord}_\ell(d) > 0$ and $\mathrm{ord}_\ell(d') > 0$. Let $\tau$ and $\tau'$ denote the images in $\Phi$ of $E(C, C')$ and $E(D, D')$. Then Proposition 3.2 implies that

$$\langle \tau, \tau' \rangle = \mathrm{lcm}(r, r')\mathrm{lcm}(d, d')(1/dd').$$

Since $\ell \nmid rr'$, $\langle \tau, \tau' \rangle$ is not trivial in $\mathbb{Q}_\ell/\mathbb{Z}_\ell$. Thus, the $\ell$-parts of $\tau$ and $\tau'$ are not trivial in $\Phi$. Note that it follows from Lemma 2.2 that $\mathrm{ord}_\ell(\tau') = \mathrm{ord}_\ell(\gcd(d, d'))$.

**Proposition 3.4** *Let $(G, M, R)$ be any arithmetical graph. Let $(C, r)$, $(C_1, r_1)$, $\ldots$, $(C_n, r_n)$ be a terminal chain $T$ of $G$, with node $C$ and terminal vertex $C_n$. Then the image $\tau$ of $E(C, C_j)$ in $\Phi(G)$ is trivial for all $j = 1, \ldots, n$.*

*Proof:* Since the pairing $\langle\ ;\ \rangle$ is perfect it is sufficient, to show that $\tau = 0$, to show that $\langle \tau, \sigma \rangle = 0$ for all $\sigma \in \Phi(G)$. Let $\sigma$ denote the image of $\Phi$ of $E(D, D')$, where $D, D'$ are any vertices of $G$, of multiplicity $r_D$ and $r_{D'}$. If neither $D$ nor $D'$ belong to the terminal chain $T$, or if $D = C$ and $D' \notin T$, then Proposition 3.2 implies that $\langle \tau, \sigma \rangle = 0$. Assume now that $D = C_i$ and $D' \neq C_s$, for all $s = 1, \ldots, n$. Let $m = \max(i, j)$. Then, using 3.2 and 2.8, we find that there exist two integers $b$ and $c$ such that

$$\begin{aligned} \langle \tau, \sigma \rangle &= \mathrm{lcm}(r, r_j)\ \mathrm{lcm}(r_i, r_{D'}) \left( \frac{1}{r_m r_{m-1}} + \ldots + \frac{1}{r_1 r} \right) \\ &= \mathrm{lcm}(r, r_j)\ \mathrm{lcm}(r_i, r_{D'}) \left( \frac{b}{r r_n} - \frac{c}{r_m r_{m-1}} \right). \end{aligned}$$

Since $r_n \mid r_i$, we have $\mathrm{lcm}(r, r_j)\ \mathrm{lcm}(r_i, r_{D'}) b / r r_n = 0$ in $\mathbb{Q}/\mathbb{Z}$. If $m = i$, then we use the fact that $r_n \mid r$ to find that $\mathrm{lcm}(r, r_j)\ \mathrm{lcm}(r_i, r_{D'}) c / r_i r_n = 0$ in $\mathbb{Q}/\mathbb{Z}$. If $m = j$, we use again the fact that $r_n \mid r_i$ to find that $\mathrm{lcm}(r, r_j)\ \mathrm{lcm}(r_i, r_{D'}) c / r_j r_n = 0$ in $\mathbb{Q}/\mathbb{Z}$. Thus, in all cases, $\langle \tau, \sigma \rangle = 0$. If $D = C_i$ and $D' = C_s$, Proposition 2.7 shows that $\sigma = 0$. This concludes the proof of Proposition 3.4. The reader may use the techniques developped in the above proof to give a different proof of Proposition 2.7.

**Remark 3.5** Let $(C_1, r_1)$, $(C_2, r_2)$, and $(D, r)$, be three vertices on an arithmetical graph $(G, M, R)$. Then

$$rE(C_1, C_2) = \frac{r_2 \gcd(r, r_1)}{\gcd(r_2, r_1)} E(C_1, D) + \frac{r_1 \gcd(r, r_2)}{\gcd(r_2, r_1)} E(D, C_2).$$

If $\ell \nmid rr_1r_2$, we find that the order of the $\ell$-part of $E(C_1, C_2)$ divides the maximum of the orders of the $\ell$-parts of $E(C_1, D)$ and $E(D, C_2)$.

If $C_1$ and $C_2$ belong to the same terminal chain of $G$ and if $\ell \nmid r_1r_2$, we find, using 2.7 and 3.4, that the orders of the $\ell$-parts of $E(C_1, D)$ and $E(C_2, D)$ are equal.

**3.6** If $D$ is a node and $C_1$ and $C_2$ are terminal vertices on two terminal chains attached to $D$, then we shall call $(C_1, C_2)$ an *elementary pair*. In the case of an elementary pair, both $r_1$ and $r_2$ divide $r$ and we find that as vectors in $\mathbb{Z}^v$,

$$\frac{r}{\operatorname{lcm}(r_1, r_2)} E(C_1, C_2) = E(C_1, D) + E(D, C_2).$$

Using Proposition 3.4, we see that $E(C_1, C_2)$ has order dividing $r/\operatorname{lcm}(r_1, r_2)$. We shall compute in 4.3 the precise order of such a pair of vertices.

# 4 Elementary pairs

**4.1** Let $(G, M, R)$ be an arithmetical graph with $v$ vertices. Let us recall how one may compute the order of the group $\Phi(G)$. Let $\tilde{R} := \operatorname{diag}(r_1, \ldots, r_v)$. Consider the arithmetical graph $(\tilde{G}, \tilde{M}, J)$, where ${}^tJ := (1, \ldots, 1)$ and $\tilde{M} := \tilde{R}M\tilde{R}$. Then (see [Lor1], 1.11)

$$|\Phi(\tilde{G})| = (r_1 \cdots r_v)^2 |\Phi(G)|.$$

This formula is useful in practice because it reduces the computation of $|\Phi(G)|$ to counting the number of spanning trees of the graph $\tilde{G}$.

Fix a numbering of the vertices of $G$. Let $e_1, \ldots, e_v$ denote the standard basis of $\mathbb{Z}^v$. Let $E_{ij} := e_i - e_j$. When two vertices of $G$ are denoted $(C, r)$ and $(C', r')$ without specifying a numbering $i$ for $C$ and $j$ for $C'$, we may use $E_{CC'}$ to denote the vector $E_{ij}$. Suppose given a vector ${}^tS = (s_1, \ldots, s_v) \in \mathbb{Z}^v$ such that

$$(\tilde{R}M\tilde{R})S = \mu E_{CC'}.$$

Then, by definition, the order of $E_{CC'}$ in $\Phi(\tilde{G})$ divides $\mu$. Note also that $r$ and $r'$ divide $\mu$. Recall that

$$ {}^tE(C,C') := (0,\ldots,0,r'/\gcd(r,r'),0,\ldots,0,-r/\gcd(r,r'),0,\ldots,0). $$

It follows that

$$ M(\tilde{R}S) = \frac{\mu}{\mathrm{lcm}(r,r')} E(C,C'). $$

Let $g$ be the greatest common divisor of the coefficients of the vector $\tilde{R}S$. Then the order of $E(C,C')$ in $\Phi(G)$ divides $\mu/g\mathrm{lcm}(r,r')$.

**Lemma 4.2** *Suppose that in $S$, the coefficient $s_k$ corresponding to the vertex $(C_k,r_k)$ is null for some $k \in \{1,\ldots,v\}$. If $\ell$ is prime and $\ell \nmid r_k$, then the $\ell$-part of the order of $E(C,C')$ is equal to the $\ell$-part of $\mu/g\mathrm{lcm}(r,r')$.*

*Proof:* Assume that the order of $E(C,C')$ in $\Phi(G)$ is equal to $d$, with $cd = \mu/g\mathrm{lcm}(r,r')$, $c \in \mathbb{N}$. Then there exists an integer vector ${}^tZ = (z_1,\ldots,z_v)$ such that $MZ = dE(C,C')$. Hence, $(\tilde{R}S/g) - cZ = fR$ for some integer $f$. From the relation $cz_k = fr_k$, we conclude that $\mathrm{ord}_\ell(c) \le \mathrm{ord}_\ell(f)$. Since $cZ + fR = \tilde{R}S/g$, we find that every coefficient of $(\tilde{R}S/g)$ is divisible by $\ell^{\mathrm{ord}_\ell(c)}$. Thus we have $\mathrm{ord}_\ell(c) = 0$, which implies that $\mathrm{ord}_\ell(d) = \mathrm{ord}_\ell(\mu/g\mathrm{lcm}(r,r'))$, as desired.

Let us now apply the above technique to compute the order of the following elementary pair $(C_n, C'_{n'})$. Let $(D,r)$ be a node of the graph $G$. Let $(D,r)$, $(C_1,r_1),\ldots,(C_n,r_n)$ be a terminal chain $T$ on $G$ with terminal vertex $C_n$. Let $(D,r)$, $(C'_1,r'_1),\ldots,(C'_{n'},r'_{n'})$ be a terminal chain $T'$ on $G$ with terminal vertex $C'_{n'}$. Let $G_D$ denote the connected component of $D$ in $G \setminus \{\text{edges of } T \cup T'\}$. Let

$$ m := \gcd(r_j, (D_j,r_j) \text{ any vertex of } G_D). $$

Note that $m \mid r$.

**Proposition 4.3** *Keep the notation introduced in this section.*

- *If $\ell \nmid r_n r'_{n'}$, then the $\ell$-part of the element $E(C_n, C'_{n'})$ has order $\ell^{\mathrm{ord}_\ell(r/m)}$ in $\Phi(G)$.*

- *If exactly one multiplicity $r_n$ or $r'_{n'}$ is divisible by $\ell$, then the $\ell$-part of the element $E(C_n, C'_{n'})$ has order $\ell^{\mathrm{ord}_\ell(r/\mathrm{lcm}(r_n, r'_{n'}))}$ in $\Phi(G)$.*

*Proof:* If $\ell \nmid r$, then the proposition follows from 3.6. Thus we assume from now on that $\ell \mid r$. Without loss of generality, we may assume that $\ell \nmid r_n$. Consider the graph $\tilde{G}$ associated to $G$. Let us explicitly write down a vector $S$ associated to the pair $(C_n, C'_{n'})$ in $\tilde{G}$. Set

$$\mu := \mathrm{lcm}(r_n r_{n-1}, r_{n-1} r_{n-2}, \ldots, r_1 r, r r'_1, \ldots, r'_{n'-1} r'_{n'}).$$

The reader will check that the following vector ${}^tS := (s_{C_n}, s_{C_{n-1}}, \ldots)$ is such that $\tilde{M}S = \mu E_{C_n C'_{n'}}$, where

$$\begin{array}{rcl} s_{C_n} & := & 0, \\ s_{C_{n-1}} & := & \mu/r_n r_{n-1}, \\ s_{C_{n-2}} & := & \mu/r_n r_{n-1} + \mu/r_{n-1} r_{n-2}, \\ \vdots & & \\ s_D & := & \mu/r_n r_{n-1} + \mu/r_{n-1} r_{n-2} + \ldots + \mu/r_1 r, \\ s_{C'_1} & := & s_D + \mu/r r'_1, \\ \vdots & & \\ s_{C'_{n'}} & := & s_D + \mu/r r'_1 + \ldots + \mu/r'_{n'-1} r'_{n'}. \\ s_C & := & s_D, \text{ if } C \text{ is any vertex of } G_D. \end{array}$$

Let $\{b_1, \ldots, b_n\}$ denote the sequence of integers associated in 2.3 to the terminal chain $T$. Lemma 2.8 shows that

$$s_D := \mu(1/r_n r_{n-1} + 1/r_{n-1} r_{n-2} + \ldots + 1/r_1 r) = \mu b_n / r r_n.$$

Lemma 2.6 shows that $b_n$ and $r/r_n$ are coprime. The same argument shows that $s_{C_i}$ equals $\mu f_i / r_i r_n$ for some integer $f_i$ with $\gcd(f_i, r_i/r_n) = 1$.

When $\ell \mid r'_{n'}$, we claim that $\ell \nmid m$. Indeed, we know from $MR = 0$ that $|(D \cdot D)| r = r_1 + r'_1 + z$, with $z$ divisible by $\ell^{\mathrm{ord}_\ell(m)}$. Since $\ell^{\mathrm{ord}_\ell(m)} \mid r$, we find that $\ell^{\mathrm{ord}_\ell(m)} \mid r_1 + r'_1$. Since $r_n = \gcd(r, r_1)$ and $r'_n = \gcd(r, r'_1)$, we find that when only one of the terminal multiplicities is divisible by $\ell$, $\mathrm{ord}_\ell(m) = 0$. Since $s_{C_n} = 0$ and $\ell \nmid r_n$, it follows from 4.2 and the above discussion that to prove both parts of the proposition, it is sufficient to show that the $\ell$-part of the greatest common divisor $g$ of the coefficients of the vector ${}^t(\tilde{R}S) = (r_n s_{C_n}, r_{n-1} s_{C_{n-1}}, \ldots)$ is equal to the $\ell$-part of $\mu m / r$.

Let $(A, a)$ be a vertex of $G_D$ with $\mathrm{ord}_\ell(a) = \mathrm{ord}_\ell(m)$. Then $\mathrm{ord}_\ell(as_A) = \mathrm{ord}_\ell(\mu b_n a/r)$. Since $\ell \nmid b_n$ (because $\ell \mid r$ and $\gcd(r/r_n, b_n) = 1$), $\mathrm{ord}_\ell(as_A) = \mathrm{ord}_\ell(\mu a/r)$. To conclude the proof of the proposition, we only need to show that every coefficient of $\tilde{R}S$ is divisible by the $\ell$-part of $\mu a/r$; i.e., that $\mathrm{ord}_\ell(g) = \mathrm{ord}_\ell(\mu a/r)$. (Note that the $\ell$-part of $\mu a/r$ divides the $\ell$-part of $\mu$.) For all $i = 1, \ldots, n$, the coefficient $r_i s_{C_i}$ equals $r_i \mu f_i / r_i r_n$, and thus is divisible by the $\ell$-part of $\mu$. The coefficient $rs_D$ equals $r\mu b_n / r r_n$ and is also divisible by the $\ell$-part of $\mu$. Recall that

$$s_{C_1'} = \frac{\mu b_n}{r r_n} + \frac{\mu}{r r_1'} = \frac{\mu f_1}{r_1 r_n} + \mu(\frac{1}{r_1 r} + \frac{1}{r r_1'}) = \frac{\mu f_1}{r_1 r_n} + \frac{\mu}{r}(\frac{r_1 + r_1'}{r_1 r_1'}).$$

We have seen above that $\ell^{\mathrm{ord}_\ell(m)} \mid r_1 + r_1'$. Since $r_n = \gcd(r, r_1)$ is prime to $\ell$, we conclude that $\ell \nmid r_1$. Hence, $r_1' s_{C_1'}$ is divisible by the $\ell$-part of $\mu a/r$. Similarly, when $\ell \nmid r_{n'}'$, then $\ell \nmid r_1'$, and in this case we find that $s_{C_1'}$ is divisible by the $\ell$-part of $\mu a/r$.

Let us conclude the proof in the case where $\ell \nmid r_{n'}'$. The coefficient $r_2' s_{C_2'} = r_2'(s_{C_1'} + \mu/r_1' r_2')$ is clearly divisible by the $\ell$-part of $\mu a/r$. Consider now

$$r_3' s_{C_3'} = r_3'(s_{C_2'} + \mu/r_2' r_3').$$

If $\ell \nmid r_2'$, then it is clear that $r_3' s_{C_3'}$ is divisible by the $\ell$-part of $\mu a/r$. If $\ell \mid r_2'$, recall that $r_2'$ divides $r_1' + r_3'$. In particular, $r_3'$ is not divisible by $\ell$. Since

$$r_3' s_{C_3'} = r_3'(s_{C_1'} + \frac{\mu}{r_2'} \frac{r_1' + r_3'}{r_1' r_3'}),$$

we are able again to conclude that $r_3' s_{C_3'}$ is divisible by the $\ell$-part of $\mu a/r$. We leave it to the reader to check that this process can be continued to show that all coefficients $r_i' s_{C_i'}$ are divisible by the $\ell$-part of $\mu a/r$.

Finally, let us discuss the case where $\ell \mid r_{n'}'$ (recall that in this case, $\mathrm{ord}_\ell(a) = 0$). The case of $r_2' s_{C_2'}$ is a variation on the case of $r_3' s_{C_3'}$, treated below. The coefficient $r_3' s_{C_3'} = r_3'(s_{C_2'} + \mu/r_2' r_3')$ is divisible by the $\ell$-part of $\mu/r$ if $\mathrm{ord}_\ell(r) \geq \mathrm{ord}_\ell(r_2')$ and $\mathrm{ord}_\ell(r_3') \geq \mathrm{ord}_\ell(r_2')$. It follows from $MR = 0$ that $r_1'$ divides $r + r_2'$. Thus, if $\mathrm{ord}_\ell(r) < \mathrm{ord}_\ell(r_2')$, we find that $\mathrm{ord}_\ell(r) \geq \mathrm{ord}_\ell(r_1')$. Since $r_2'$ divides $r_1' + r_3'$, we find that $\mathrm{ord}_\ell(r_3') = \mathrm{ord}_\ell(r_1')$.

If $\mathrm{ord}_\ell(r) \geq \mathrm{ord}_\ell(r_2')$ and $\mathrm{ord}_\ell(r_3') < \mathrm{ord}_\ell(r_2')$, then $\mathrm{ord}_\ell(r_1') = \mathrm{ord}_\ell(r_3')$. These conditions imply that $\mathrm{ord}_\ell(r) \geq \mathrm{ord}_\ell(r_1')$ and $\mathrm{ord}_\ell(r_3') \geq \mathrm{ord}_\ell(r_1')$.

Write now

$$s_{C'_3} = s_{C'_1} + \frac{\mu}{r'_2}\frac{r'_3 + r'_1}{r'_1 r'_3}.$$

It is easy to see that the coefficient $r'_3 s_{C'_3}$ is divisible by the $\ell$-part of $\mu/r$ if $\text{ord}_\ell(r) \geq \text{ord}_\ell(r'_1)$ and $\text{ord}_\ell(r'_3) \geq \text{ord}_\ell(r'_1)$. Hence, we conclude that $r'_3 s_{C'_3}$ is divisible by the $\ell$-part of $\mu/r$. We leave it to the reader to check that this process can be continued to show that all coefficients $r'_i s_{C'_i}$ are divisible by the $\ell$-part of $\mu/r$. This concludes the proof of Proposition 4.3.

**Example 4.4** Consider the following graphs.

The pair $E(C, C')$ has order 4 in $\Phi(G_1)$, and order 2 in $\Phi(G_2)$.

Using the notation introduced in Proposition 4.3, let $\tau \in \Phi(G)$ denote the image of $E(C_n, C'_{n'})$. It follows from Proposition 3.2 that

$$\begin{aligned} \langle \tau, \tau \rangle &= \text{lcm}(r_n, r'_{n'})^2 (1/r_n r_{n-1} + \ldots + 1/r_1 r + 1/r r'_1 + \ldots + 1/r'_{n'-1} r'_{n'}) \\ &= \text{lcm}(r_n, r'_{n'})^2 (b_n / r r_n + b'_{n'} / r r'_{n'}). \end{aligned}$$

Recall that in $G$, we have $r_1 + r'_1 + z = \Delta r$, for some $\Delta \in \mathbb{N}$.

**Proposition 4.5** *Keep the notation introduced above.*

- *If $\ell \nmid r_n r'_{n'}$, then the $\ell$-part of the order of $\langle \tau, \tau \rangle$ in $\mathbb{Q}/\mathbb{Z}$ is equal to $\ell^{\text{ord}_\ell(r/z)}$. In particular this $\ell$-part may be smaller than $\ell^{\text{ord}_\ell(r/m)}$, the $\ell$-part of the order of $\tau$ in $\Phi(G)$.*

- *If exactly one multiplicity $r_n$ or $r'_{n'}$ is divisible by $\ell$, then the $\ell$-part of the order of $\langle \tau, \tau \rangle$ in $\mathbb{Q}/\mathbb{Z}$ is equal to $\ell^{\text{ord}_\ell(r/\text{lcm}(r_n, r'_{n'}))}$. Thus, this $\ell$-part equals the $\ell$-part of the order of $\tau$ in $\Phi(G)$.*

*Proof:* We may assume that $\mathrm{ord}_\ell(r'_{n'}) = 1$. Using 3.6, we find that if $\mathrm{ord}_\ell(r/r_n) = 1$, then the $\ell$-part of $\tau$ is trivial, and thus Proposition 4.5 holds. Assume now that $\mathrm{ord}_\ell(r/r_n) > 1$. Recall that

$$\langle \tau, \tau \rangle = \frac{\mathrm{lcm}(r_n, r'_{n'})^2}{r_n r'_{n'}} \left( \frac{b_n r'_{n'} + b'_{n'} r_n}{r} \right).$$

Lemma 2.6 shows that $b_n r_1 - r_n = er$ and $b'_{n'} r'_1 - r'_{n'} = fr$ for some integers $e, f$. Thus

$$\begin{aligned} r'_1(b_n r'_{n'} + b'_{n'} r_n) &= r'_1 b_n r'_{n'} + (fr + r'_{n'}) r_n \\ &= (\Delta r - z - r_1) b_n r'_{n'} + (fr + r'_{n'}) r_n \\ &= (\Delta r - z) b_n r'_{n'} - (er + r_n) r'_{n'} + (fr + r'_{n'}) r_n. \end{aligned}$$

It follows that modulo $\mathbb{Z}$,

$$r'_1(b_n r'_{n'} + b'_{n'} r_n)/r \equiv -z b_n r'_{n'}/r.$$

Since $r_n = \gcd(r_1, r)$ and $b_n r_1 - 1 = ar$, we find that $\gcd(\ell, b_n) = 1$. Similarly, since $\mathrm{ord}_\ell(r'_{n'}) = 1$, $\mathrm{ord}_\ell(r'_1) = 1$. If $\mathrm{ord}_\ell(r'_{n'}) = 1$ and $\mathrm{ord}_\ell(r_n) > 1$, then the relation $r_1 + r'_1 + z = \Delta r$ shows that in this case $\mathrm{ord}_\ell(z) = 1$. Proposition 4.5 follows.

**Corollary 4.6** *If $\ell \nmid r_n r'_{n'}$ and $\mathrm{ord}_\ell(z) = \mathrm{ord}_\ell(m)$, or if exactly one multiplicity $r_n$ or $r'_{n'}$ is divisible by $\ell$, then $\tau$ is not divisible by $\ell$ in $\Phi$.*

*Proof:* Given any two elements $\tau$ and $\sigma$ of $\Phi$ of orders $t$ and $s$, respectively, it is easy to check that the order of the element $\langle \tau, \sigma \rangle$ divides $\gcd(t, s)$. Hence, if $\tau = \ell \xi$ in $\Phi$, then $\langle \tau, \tau \rangle = \ell \langle \tau, \xi \rangle$ is killed by $t/\ell$. Now let $\tau$ be as in 4.5. Since Proposition 4.5 shows that $\ell^{\mathrm{ord}_\ell(t)}$ divides the order of $\langle \tau, \tau \rangle$, we find that $\tau$ is not divisible by $\ell$.

The reader will find in 7.6 an example where $\ell \nmid r_n r'_{n'}$ and $\mathrm{ord}_\ell(z) > \mathrm{ord}_\ell(m)$, and where $\tau$ is divisible by $\ell$ in $\Phi$.

# 5 A splitting of the group of components

Let $(G, M, R)$ be an arithmetical graph. Fix a prime $\ell$. Let $\Phi_\ell(G)$ denote the $\ell$-part of the group of components $\Phi(G)$. Let $(D, r)$ be a vertex of $G$ such that $G \setminus \{D\}$ is not connected. Our aim in this section is to establish an isomorphism between $\Phi_\ell(G)$ and the product of the $\ell$-parts of the groups of components of arithmetical graphs associated to the connected components of $G \setminus \{D\}$.

**Construction 5.1** Label the connected components of $G \setminus \{D\}$ by $\mathcal{G}_1, \ldots, \mathcal{G}_t$. Label the vertices of $\mathcal{G}_i$ adjacent in $G$ to $D$ by $(C_{i,1}, r_{i,1}), \ldots, (C_{i,s_i}, r_{i,s_i})$. Assume that $t > 1$. For $i = 1, \ldots, t$, let $g_i$ denote the greatest common divisor of $r$ and the multiplicities of all vertices of $\mathcal{G}_i$.

Construct a new connected arithmetical graph $G_i$ associated to $\mathcal{G}_i$ as follows. Start with $\mathcal{G}_i \cup \{D\}$. Give to $D$ the multiplicity $r/g_i$. Give to a vertex in $\mathcal{G}_i$ its multiplicity in $G$ divided by $g_i$. Let $c_i$ denote the least integer such that $c_i r - \sum_{j=1}^{s_i} (C_{ij} \cdot D) r_{ij} \geq 0$. The integer $c_i$ will be the self-intersection of $D$ in $G_i$.

If $r$ divides $\sum_{j=1}^{s_i} (C_{ij} \cdot D) r_{ij}$, then the graph $G_i := \mathcal{G}_i \cup \{D\}$ with multiplicities as above is an arithmetical graph. If $r$ does not divide $\sum_{i=1}^{s_i} (C_{ij} \cdot D) r_{ij}$, then let $\hat{r}_i := (c_i r - \sum_{j=1}^{s_i} (C_{ij} \cdot D) r_{ij})/g_i$. Construct a terminal chain $T$ using $(r/g_i, \hat{r}_i)$ and Euclid's algorithm as in [Lor2], 2.4. The graph $G_i$ consists then in the graph $\mathcal{G}_i \cup \{D\}$, with the chain $T$ attached to $D$.

We shall say that the graph $G$ is *$\ell$-breakable at* $(D, r)$ if $\ell \nmid r$ and $t > 1$.

**Example 5.2** Let $G$ be the following graph.

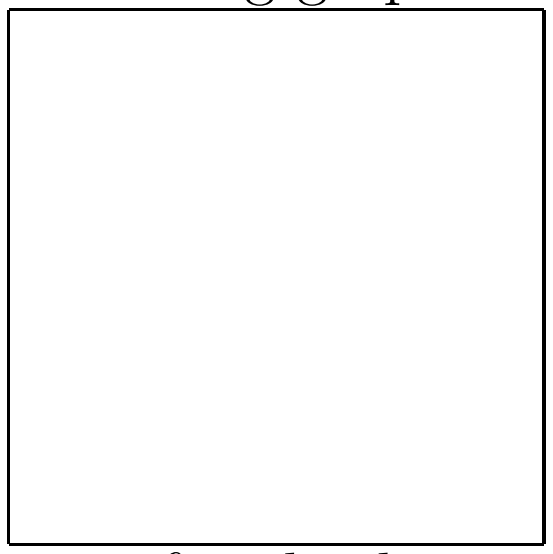

Let $D$ denote the central vertex of multiplicity 6. Then $G \setminus \{D\}$ has 3 components $\mathcal{G}_1$, $\mathcal{G}_2$ and $\mathcal{G}_3$, and the above procedure produces 3 new arithmetical graphs:

Note that $|\Phi(G_1)| = 2$, $|\Phi(G_2)| = 8$, and $|\Phi(G_3)| = 1$. Our next proposition shows then that the only primes that can divide $|\Phi(G)|$ are 2 or 3. The orders of these groups of components are best computed using the method recalled in 4.1. One finds in particular that the full group $\Phi(G)$ has order 144.

**Proposition 5.3** *Let $(G, M, R)$ be any arithmetical graph. Let $\ell$ be a prime. Assume that $G$ is $\ell$-breakable at a vertex $(D, r)$. Let $G_1, \ldots, G_t$, denote the arithmetical graphs associated as in 5.1 to the components of $G \setminus \{D\}$. Then there exists an isomorphism*

$$\alpha : \Phi_\ell(G) \longrightarrow \prod_{i=1}^{t} \Phi_\ell(G_i).$$

*Let $(C_1, r_1)$ and $(C_2, r_2)$ be any two vertices of $G$. If $C_1$ and $C_2$ belong to the same component of $G \setminus \{D\}$, say to $\mathcal{G}_j$, or if $C_1 \in \mathcal{G}_j$ and $C_2 = D$, then we denote by $E(C_1, C_2)$ and $E(C_i, D)$ both the elements of $\Phi(G)$ and the corresponding elements of $\Phi(G_j)$. Then the $\ell$-part of $E(C_1, C_2)$ is mapped under $\alpha$ to the element of $\prod_{i=1}^{k} \Phi_\ell(G_i)$ having the $\ell$-part of $E(C_1, C_2)$ in the $j$-th coordinate, and 0 everywhere else.*

*Proof:* Let $N_i$, $i = 1, \ldots, t$, denote the square submatrix of the intersection matrix $M$ corresponding to the vertices of $G$ that belong to $\mathcal{G}_i$. The matrix $M$ has the following form:

$$\begin{pmatrix} (D \cdot D) & * & \ldots & \ldots & * \\ * & N_1 & 0 & \ldots & 0 \\ \vdots & 0 & N_2 & 0 & \vdots \\ \vdots & \vdots & & \ddots & 0 \\ * & 0 & \ldots & 0 & N_t \end{pmatrix}$$

(in particular, the first column is the '$D$-column'). Multiply the first row by $r$, and add to it the sum of all other rows, each multiplied by its corresponding multiplicity. Perform a similar operation on the first column of $M$, to obtain a matrix $M'$ of the form

$$\begin{pmatrix} 0 & 0 & \dots & \dots & 0 \\ 0 & N_1 & & & \\ \vdots & & N_2 & & \\ \vdots & & & \ddots & \\ 0 & & & & N_t \end{pmatrix}.$$

Since $\ell \nmid r$, the row and column operations described above are permissible over $\mathbb{Z}_\ell$. The module $\mathrm{Ker}({}^tR) \otimes_{\mathbb{Z}} \mathbb{Z}_\ell$ is thus the direct sum of $t$ $\mathbb{Z}_\ell$-submodules $V_i$, $i = 1, \dots, t$, where

$$V_i := \oplus_{C \in \mathcal{G}_i} \mathbb{Z}_\ell r^{-1} E(C, D).$$

Let $W_i$ denote the $\mathbb{Z}_\ell$-span of the column vectors of $N_i$, so that $\mathrm{Im}(M) \otimes \mathbb{Z}_\ell \cong \oplus_{i=1}^t W_i$. Then

$$\Phi_\ell(G) \cong \oplus_{i=1}^t V_i / W_i.$$

We claim that $\Phi_\ell(G_i) \cong V_i/W_i$. Indeed, the intersection matrix $M_i$ associated to $G_i$ has the following form

$$\left(\begin{array}{ccccc|cc} * & 1 & & & & & \\ 1 & \ddots & \ddots & & & & \\ & \ddots & * & 1 & & & \\ & & 1 & (D \cdot D) & & * & \dots \\ \hline & & & * & & & \\ & & & \vdots & & & N_i \end{array}\right) \quad \text{or} \quad \begin{pmatrix} (D \cdot D) & \dots & * \\ \vdots & & \\ * & \boxed{N_i} & \end{pmatrix},$$

where the case on the right occurs if $G_i = \mathcal{G}_i \sqcup \{D\}$ (see 5.1). Multiply the $D$-row by $r/g_i$, and add to it all other rows multiplied by their corresponding multiplicities. Perform a similar operation on the columns of $M_i$ to get a matrix $M_i'$ of the form

$$\begin{pmatrix} * & 1 & & & 0 & & \\ 1 & \ddots & \ddots & & \vdots & & \\ & \ddots & \ddots & 1 & \vdots & & \\ & & 1 & * & 0 & & \\ 0 & \dots & \dots & 0 & 0 & \dots & 0 \\ & & & & \vdots & \multicolumn{2}{|c}{N_i} \\ & & & & 0 & \multicolumn{2}{|c}{} \end{pmatrix} \quad \text{or} \quad \begin{pmatrix} 0 & \dots & 0 \\ \vdots & & \\ 0 & \boxed{N_i} & \end{pmatrix}.$$

In the case of the first matrix $M_i'$, the top left corner can be further reduced to:

$$\begin{pmatrix} 0 & & & -\alpha & 0 & & \\ 1 & \ddots & & 0 & \vdots & & \\ & \ddots & \ddots & \vdots & \vdots & & \\ & & 1 & 0 & 0 & & \\ 0 & \dots & \dots & 0 & 0 & \dots & 0 \\ & & & & \vdots & \multicolumn{2}{|c}{N_i} \\ & & & & 0 & \multicolumn{2}{|c}{} \end{pmatrix},$$

where $\alpha = -r/\gcd(r, \hat{r}_i)$ (see the proof of 2.7; note that $\gcd(r, \hat{r}_i)$ is the terminal multiplicity of the terminal chain attached to $D$ in $G_i$). Since $\ell \nmid r$, we find that in both cases,

$$\Phi_\ell(G_i) \cong \oplus_{C \in \mathcal{G}_i} \mathbb{Z}_\ell r^{-1} E(C, D)/\mathrm{Im}(N_i).$$

Therefore, we have an isomorphism between $\Phi_\ell(G_i)$ and $V_i/W_i$, where the $\ell$-part of $E(C, D)$ in $\Phi_\ell(G_i)$ is mapped to the $\ell$-part of the element $E(C, D)$ in $V_i/W_i \subset \Phi_\ell(G)$. We leave it to the reader to compute the image of $E(C_1, C_2)$ under the above isomorphism. This concludes the proof of Proposition 5.3.

We now use Proposition 5.3 to prove the following important theorem.

**Theorem 5.4** *Let $(G, M, R)$ be any arithmetical graph. Let $\ell$ be any prime. Let $(C, r)$ and $(C', r')$ be a weakly connected and $\ell$-breakable pair of vertices of $G$ with $\ell \nmid rr'$ and associated integer $\lambda(C, C')$ as in 1.2. Then the $\ell$-part of $E(C, C')$ has order $\lambda(C, C')$ in $\Phi(G)$.*

*Proof:* Suppose that both $C$ and $C'$ are not terminal vertices of $G$. Then by hypothesis the graph is $\ell$-breakable at $C$ into two or more arithmetical graphs. Denote by $G'$ the new arithmetical graph that contains $C'$ (graph constructed while breaking $G$). Then the $\ell$-part of $E(C, C')$ in $G$ has the same order as the $\ell$-part of $E(C, C')$ in $G'$. Since $(C, C')$ is weakly connected, $C$ lies on a terminal chain $T$ of $G'$, and its multiplicity is still coprime to $\ell$. Moreover, $(C, C')$ is a weakly connected $\ell$-breakable pair of $G'$. Denote by $D$ the terminal vertex of $T$. If $C \neq D$, we find using Remark 3.5 that the $\ell$-part of $E(C, C')$ has the same order as the $\ell$-part of $E(D, C')$. The pair $(D, C')$ is clearly a weakly connected $\ell$-breakable pair of $G'$. Thus, to prove Theorem 5.4 for $(C, C')$ in $G$, it is sufficient to prove it for pairs where one vertex is a terminal vertex, say the vertex $C$.

Let $\mathcal{P}$ denote the path associated to the weakly connected pair $(C, C')$. (Note that if $\mathcal{P} \setminus \{C, C'\}$ contains no vertices, the theorem follows from 2.2.) Let $(C_1, r_1), (C_2, r_2), \ldots, (C_s, r_s)$ be the nodes on $\mathcal{P} \setminus \{C, C'\}$, as discussed in 1.1. If $s = 0$, then $C$ and $C'$ belong to the same terminal chain of $G$, and Theorem 5.4 follows from 2.7. If $s = 1$ and $C'$ is not a terminal vertex, we may apply the reduction step described at the beginning of the proof and assume without loss of generality that $C'$ is a terminal vertex. Then we can apply 4.3 to show that our statement holds in this case.

We proceed by induction on the number $s$ of nodes on $\mathcal{P}$. Let $m > 1$ and assume that Theorem 5.4 holds for $s \leq m - 1$. Let $(C, C')$ be a pair whose associated path $\mathcal{P}$ contains $m$ nodes. Since the pair is $\ell$-breakable, there exists a vertex $(D, r_D)$ on $\mathcal{P}$ with $\ell \nmid r_D$ and such that both components of $\mathcal{P} \setminus \{D\}$ contain at most $m - 1$ nodes $C_i$. (Note that one component may contain no nodes such as, for instance, when $r_D = C_1$.) Break the graph $G$ at $D$. Call $G_1$ the arithmetical graph associated to the connected component of $G \setminus \{D\}$ which contains $C$. Call $G_2$ the arithmetical graph that contains $C'$. The pairs $(C, D)$ and $(D, C')$ are weakly connected and $\ell$-breakable pairs of $G_1$ and $G_2$, respectively. We may thus apply the induction hypothesis to both pairs. To conclude the proof of Theorem 5.4, we need only to show that the order of the $\ell$-part of $E(C, C')$ is equal to the maximum of the orders of the $\ell$-parts of $E(C, D)$ and $E(D, C')$. Note that 3.5 only shows that the order of $E(C, C')$ divides the maximun of the orders of $E(C, D)$ and $E(D, C')$. To prove our claim, we need to use the fact that breaking the graph $G$ at $D$ produces a splitting of $\Phi(G)_\ell$, with $\Phi(G_1)_\ell$ and $\Phi(G_2)_\ell$ as direct summands.

## 6 The subgroups $\Psi_{K,L}$ and $\Theta_K^{[3]}$

Let $X/K$ be a curve. We recall below Raynaud's description of the group $\Phi_K$ and of the map $\pi$ in terms of a regular model $\mathcal{X}/\mathcal{O}_K$ of $X/K$. Let $\mathcal{X}_k = \sum_{i=1}^v r_i C_i$, and assume that $\gcd(r_1, \ldots, r_v) = 1$. Let $\mathcal{L} := \oplus_{i=1}^v \mathbb{Z} C_i$ denote the free abelian group generated by the components $C_i$, $i = 1, \ldots, v$. Let $\mathcal{L}^* := \mathrm{Hom}_{\mathbb{Z}}(\mathcal{L}, \mathbb{Z})$, and let $\{x_1, \ldots, x_v\}$ denote the dual basis of $\mathcal{L}$, so that $x_i(C_j) = \delta_{ij}$. Let ${}^tR : \mathcal{L}^* \to \mathbb{Z}$ be the map $\sum_{i=1}^v a_i x_i \mapsto \sum_{i=1}^v a_i r_i$. Consider the following diagram.

$$\begin{array}{ccccccc}
\mathcal{L} & \xrightarrow{i} & \mathrm{Pic}(\mathcal{X}) & \xrightarrow{res} & \mathrm{Pic}(X) & \xrightarrow{deg} & \mathbb{Z} \\
|| & & \downarrow \phi & & \downarrow \psi & & || \\
\mathcal{L} & \xrightarrow{\mu} & \mathcal{L}^* & \longrightarrow & \mathcal{L}^*/\mu(\mathcal{L}) & \xrightarrow{{}^tR} & \mathbb{Z}
\end{array}$$

The map $i$ is defined as follows: $i(C_j) :=$ curve $C_j$ in $\mathcal{X}$, where the curve $C_j$ is viewed as an element of $\mathrm{Pic}(\mathcal{X})$. The map $res$ restricts a divisor of $\mathcal{X}$ to the open set $X$ of $\mathcal{X}$. The map $res$ is surjective because the scheme $\mathcal{X}$ is regular. The map $deg$ is defined as follows: $\deg(\sum_{i=1}^s a_i P_i) := \sum_{i=1}^s a_i [K(P_i) : K]$, where $K(P_i)$ is the residue field of $P_i$ in $X$. We denote by $\mathrm{Pic}^0(X)$ the kernel of the map $deg$. The intersection matrix $M$ of $\mathcal{X}_k$ can be thought of as a bilinear map on $\mathcal{L} \times \mathcal{L}$ and, therefore, induces a map $\mu : \mathcal{L} \to \mathcal{L}^*$ defined by $\mu(C_i) := \sum_{j=1}^v (C_i \cdot C_j) x_j$. Then ${}^tR \circ \mu = 0$. Let $D$ be an irreducible divisor on $\mathcal{X}$, and define $\phi(D) := \sum_{j=1}^v (C_j \cdot D) x_j$. The map $\psi$ is the natural map induced by $\phi$. It is well-known that the diagram above is commutative.

One easily checks that $\mathrm{Ker}({}^tR)/\mu(\mathcal{L})$ is the torsion subgroup of $\mathcal{L}^*/\mu(\mathcal{L})$. Raynaud [BLR], 9.6, showed that the group of components $\Phi_K$ of the jacobian $A/K$ of the curve $X/K$ is isomorphic to the group $\mathrm{Ker}({}^tR)/\mu(\mathcal{L})$. It follows from this description that the group $\Phi_K$ can be explicitly computed using a row and column reduction of the intersection matrix $M$ (see [Lor1], 1.4). Since the residue field $k$ is algebraically closed, $A(K) = \mathrm{Pic}^0(X)$. Raynaud ([BLR], 9.5/9 and 9.6/1) has shown that the reduction map $\pi : A(K) \to \Phi$ corresponds to the restricted map $\psi : \mathrm{Pic}^0(X) \to \mathrm{Ker}({}^tR)/\mu(\mathcal{L})$. Thus, given two points $P$ and $Q$ in $X(K)$, the image of $P - Q$ in the group $\Phi_K$ can be identified with the image of $E(C_P, C_Q)$ in $\Phi(G)$.

To prove Theorem 6.5 below, we need to recall the following facts.

**6.1** Let $\mathcal{X}/\mathcal{O}_K$ be any regular model of $X/K$. Let $M/K$ be a finite extension. Let $\mathcal{Y}/\mathcal{O}_M$ denote the normalization of the scheme $\mathcal{X} \times_{\mathrm{Spec}\,(\mathcal{O}_K)} \mathrm{Spec}\,(\mathcal{O}_M)$. Let $b : \mathcal{Y} \to \mathcal{X}$ denote the composition of the natural maps

$$\mathcal{Y} \to \mathcal{X} \times_{\mathrm{Spec}\,(\mathcal{O}_K)} \mathrm{Spec}\,(\mathcal{O}_M) \to \mathcal{X}.$$

Let $\rho : \mathcal{Z} \to \mathcal{Y}$ denote the minimal desingularization of $\mathcal{Y}$. To recall the descriptions of the maps $\rho$ and $b$, we need the following definition. Let $C_1, ..., C_m$ be irreducible components of the special fiber $\mathcal{X}_k$. The divisor $C := \sum_{i=1}^{m} C_i$ is said to be a *Hirzebruch-Young string* if the following four conditions hold: 1) $g(C_i) = 0$, for all $i = 1, ..., m$, and 2) $(C_i \cdot C_i) \leq -2$, for all $i = 1, ..., m$, and 3) $(C_i \cdot C_j) = 1$ if $|i - j| = 1$, and 4) $(C_i \cdot C_j) = 0$ if $|i - j| > 1$.

Recall that two curves $C$ and $C'$ of $\mathcal{X}$ meet at a point $P$ with normal crossings if the local intersection number $(C, C')_P$ is equal to 1. In particular, $P$ is a smooth point on both $C$ and $C'$. We say that two effective divisors meet with normal crossings if they meet with normal crossings at each intersection point.

Given any integer $m$ prime to $p$, let us denote by $M_m/K$ the unique Galois extension of $K$ of degree $m$. We shall call $M/K$ an $\ell$-extension of $K$ if $[M : K]$ is a power of $\ell$.

The following facts are well known; we state them without proof (see for instance [BPV], Theorem 5.2, when $\mathcal{X}/\mathbb{C}$ is a surface.)

**Facts 6.2** Let $q$ be a prime, $q \neq p$. Let $M := M_q$. Let $(C, r)$ be a component of $\mathcal{X}_k$.

- The map $b : \mathcal{Y} \to \mathcal{X}$ is ramified over the divisor $R := \sum_{\gcd(q, r_i) = 1} C_i$.
- Let $P \in \mathcal{Y}$ be a point such that $b(P) \in C \subset R$. If $b(P)$ is a smooth point of $R$, then $P$ is regular on $\mathcal{Y}$.
- Let $P \in \mathcal{Y}$ be such that $b(P)$ is the intersection point of two components $C$ and $C'$ of $R$. If $C$ and $C'$ meet with normal crossings at $P$, then the divisor $\rho^{-1}(P) := \sum_{i=1}^{m(P)} E_i$ is a Hirzebruch-Young string. Let $P \in D \cap D'$, where $D$ and $D'$ are irreducible components of $\mathcal{Y}_k$. Write $\tilde{D}$ for the strict transform of $D$ in $\mathcal{Z}$. Then:

$$(\rho^{-1}(P) \cdot \tilde{D}) = (E_1 \cdot \tilde{D}) = 1 = (E_{m(P)} \cdot \tilde{D}') = (\rho^{-1}(P) \cdot \tilde{D}').$$

Moreover, $(\rho^{-1}(P) \cdot E) = 0$ if $E$ is an irreducible component of $\mathcal{Z}_k$ with $E \neq \tilde{D}, \tilde{D}'$.

- If $q \nmid r$, then $b^{-1}(C) =: D$ is irreducible and the restricted map $b_{|D} : D \to C$ is an isomorphism. The curve $D$ has multiplicity $r$ in $\mathcal{Y}_k$.

- If $q \mid r$ and $C \cap R \neq \emptyset$, then $b^{-1}(C) =: D$ is irreducible and the restricted map $b_{|D} : D \to C$ is a morphism of degree $q$ ramified over $|C \cap R|$ points of $C$. The curve $D$ has multiplicity $r/q$ in $\mathcal{Y}_k$. When $C$ is smooth and meets $R$ with normal crossings, then $D$ is smooth and its genus of $D$ is computed using the Riemann-Hurwitz formula.

- If $q \mid r$ and $C \cap R = \emptyset$, then $b : b^{-1}(C) \to C$ is an etale map and each irreducible component of $b^{-1}(C)$ has multiplicity $r/q$ in $\mathcal{Y}_k$. If $b^{-1}(C)$ is not irreducible, then it is equal to the disjoint union $D_1 \sqcup \ldots \sqcup D_q$ of $q$ irreducible curves, and each restricted map $b_{|D_j} : D_j \to C$ is an isomorphism.

**Lemma 6.3** *Let $X/K$ be a curve with a regular model $\mathcal{X}/\mathcal{O}_K$ and associated arithmetical graph $(G, M, R)$. Let $\ell \neq p$ be a prime. Let $(C, C')$ be a weakly connected pair with $\ell \nmid rr'$. Let $M := M_\ell$ and consider the associated map $b \circ \rho : \mathcal{Z} \to \mathcal{X}$. Denote again by $C$ and $C'$ the strict transforms of $C$ and $C'$ in $\mathcal{Z}$.*

- *Assume that the pair $(C, C')$ is $\ell$-breakable in $G$. Then $(C, C')$ is also weakly connected and $\ell$-breakable in the graph of $\mathcal{Z}$.*

- *Assume that the pair $(C, C')$ is not $\ell$-breakable. Then $(C, C')$ is a multiply connected pair of $\mathcal{Z}$.*

*Proof:* Note that, by definition of weakly connected, two curves of the path $\mathcal{P}$ between $C$ and $C'$ that interesect do intersect with normal crossings. (Note on the other hand that our hypothesis allows other singularities on each components.) The lemma follows immediately from 6.2.

**Lemma 6.4** *Let $A/K$ be an abelian variety. Let $\tau \in \Phi_K$. Then $\tau \in \Psi_{K,L}$ if and only if there exists a finite extension $M/K$ such that $\tau \in \Psi_{K,M}$.*

*Proof:* It is clear that $\Psi_{K,M} \subseteq \Psi_{K,ML}$. It follows from the fact that $A_L/L$ has semistable reduction that the canonical map $\Phi_L \to \Phi_{ML}$ is injective. Hence, $\Psi_{K,M} \subseteq \Psi_{K,L}$.

**Theorem 6.5** *Let $X/K$ be a curve with a regular model $\mathcal{X}/\mathcal{O}_K$ and associated arithmetical graph $(G, M, R)$. Let $\ell \neq p$ be a prime. Let $(C, C')$ be a weakly connected $\ell$-breakable pair with $\ell \nmid rr'$. Then the $\ell$-part of the image of $E(C, C')$ in $\Phi_K(\mathrm{Jac}(X))$ belongs to $\Psi_{K,L}$ and has order $\lambda(C, C')$.*

*Proof:* Let $\mathcal{P}$ denote the path linking $C$ and $C'$. Let $M/K$ be any $\ell$-extension. Let $b \circ \rho : \mathcal{Z} \to \mathcal{X}$ be the associated base change and desingularization map as in 6.1. Denote again by $C$ and $C'$ the preimages in $\mathcal{Z}$ of the components $C$ and $C'$ in $\mathcal{X}$. It follows from 6.3 that the pair $(C, C')$ is also weakly connected and $\ell$-breakable in $\mathcal{Z}$. Thus, its order can be computed using Theorem 5.4.

Let us now consider the nodes on the path $\mathcal{P}'$ linking $C$ and $C'$ in $\mathcal{Z}$. Clearly, if $D$ is a vertex of $\mathcal{P}'$ such that $(b \circ \rho)(D)$ is a node of $\mathcal{P}$, then $D$ is a node on $\mathcal{P}'$. Moreover, if $D$ is a node on $\mathcal{P}'$ such that $(b \circ \rho)(D)$ is not a node of $\mathcal{P}$, then the component $(b \circ \rho)(D)$ is not smooth. The reader will note that after an extension of degree $\ell^d$, the multiplicity of the preimage in $\mathcal{Z}$ of a component $(C, r)$ on the path $\mathcal{P}$ is equal to $r\ell^{-\min(d, \mathrm{ord}_\ell(r))}$.

Define $\mu$ to be the power of $\ell$ such that

$$\mathrm{ord}_\ell(\mu) := \max\{\mathrm{ord}_\ell(r), (C, r) \text{ a component on } \mathcal{P}\}.$$

It follows from the above discussion and from Theorem 5.4 that over $M_\mu$, the pair $(C, C')$ has trivial $\ell$-part in $\Phi_{M_\mu}$. Therefore, the $\ell$-part of $E(C, C')$ in $\Phi_K$ belongs to $\Psi_{K,M_\mu}$. Thus, Lemma 6.4 implies that the $\ell$-part of the image of $E(C, C')$ in $\Phi_K(\mathrm{Jac}(X))$ belongs to $\Psi_{K,L}$. Note that it is not always true that $M_\mu \subseteq L$. This concludes the proof of Theorem 6.5.

**6.6** Let us recall now the description of the first functorial subgroup of $\Phi_{K,\ell}$ appearing in the filtration

$$\Theta_{K,\ell}^{[3]} \subseteq \Psi_{K,L,\ell} \subseteq \Theta_{K,\ell} \subseteq \Phi_{K,\ell}$$

introduced in [Lor3], 3.21. Let $A/K$ be an abelian variety. Let $T_\ell$ denote the Tate module $T_\ell A$, $\ell \neq p$. Let $\mathbb{D}_\ell := \mathbb{Q}_\ell/\mathbb{Z}_\ell$. Let $I_K := I(\overline{K}/K)$. There is a

natural isomorphism

$$\phi_{K,\ell} : \Phi_{K,\ell} \longrightarrow E := \frac{(T_\ell \otimes \mathbb{D}_\ell)^{I_K}}{T_\ell^{I_K} \otimes \mathbb{D}_\ell}.$$

Given any submodule $X$ of $T_\ell$, let $f_X : X \otimes \mathbb{Q}_\ell \to T_\ell \otimes \mathbb{D}_\ell$ denote the natural map. We denote by $t(X)$ the subgroup of $E$ generated by the elements $x \in (T_\ell \otimes \mathbb{D}_\ell)^{I_K}$ such that there exists $\tilde{x} \in X \otimes \mathbb{Q}_\ell$ with $f_X(\tilde{x}) = x$. Consider the submodules $W_{\ell,L} \subseteq T_\ell^{I_L} \subseteq T_\ell$, where $W_{\ell,L}$ is canonically isomorphic to the Tate module of the maximal torus $\mathcal{T}_L$ in the connected component of the Néron model of $A_L/L$. Then, by definition ([Lor3], 3.8),

$$\phi_{K,\ell}(\Theta^{[3]}_{K,\ell}) = t(W_{\ell,L}) \text{ and } \phi_{K,\ell}(\Psi_{K,L,\ell}) = t(T_\ell^{I_L}).$$

As we shall see in 6.10, the description of the elements of $\Theta^{[3]}_{K,\ell}$ seems to be more complicated than the description of the elements of $\Psi_{K,L,\ell}$.

Let $M/K$ be any finite separable extension. Denote by $\mathcal{A}_M/\mathcal{O}_M$ the Néron model of $A_M/M$. Let $\mathcal{A}_{M,k}/k$ denote its special fiber, with connected component $\mathcal{A}^0_{M,k}$. Let $\mathcal{T}_M \subset \mathcal{A}^0_{M,k}$ denote the maximal torus of $\mathcal{A}^0_{M,k}$. Denote by $\pi_M : A(M) \to \mathcal{A}_{M,k}(k)$ the reduction map.

**Lemma 6.7** *Let $A/K$ be an abelian variety with purely additive reduction. Let $\ell \neq p$ be any prime. Let $\tau \in \Phi_{K,\ell}$. Let $t$ denote the unique element of $A(K)_{\mathrm{tors},\ell}$ whose image in $\Phi_K$ is $\tau$. The element $\tau$ belongs to the subgroup $\Theta^{[3]}_{K,\ell}$ if and only if there exists a finite separable extension $M/K$ such that $\pi_M(t)$ belongs to $\mathcal{T}_M$.*

*Proof:* Let us first note that our hypothesis implies that $\pi_L(t)$ belongs to $\mathcal{T}_L$. Indeed, it follows from the properties of smooth connected commutative groups that the natural map $\mathcal{A}_M \to \mathcal{A}_{ML}$ restricts to a map $\mathcal{T}_M \to \mathcal{T}_{ML}$. Thus $\pi_{ML}(t)$ belongs to $\mathcal{T}_{ML}$. In particular, the image of $\tau$ under the natural map $\Phi_K \to \Phi_{ML}$ is trivial. Since the map $\Phi_L \to \Phi_{ML}$ is injective, we conclude that $\pi_L(t) \in \mathcal{A}^0_{L,k}(k)$. Since $\mathcal{A}^0_{L,k} = \mathcal{A}^0_{ML,k}$ by semistability, we find that $\pi_L(t) \in \mathcal{T}_L$.

When $A$ has purely additive reduction, $T_\ell^{I_K} = (0)$. The canonical reduction map $A(K)_{\mathrm{tors},\ell} \to \Phi_{K,\ell}$ is an isomorphism and factors through $(T_\ell \otimes \mathbb{D}_\ell)^{I_K}$ as follows:

$$A(K)_{\mathrm{tors},\ell} \xrightarrow{g} (T_\ell \otimes \mathbb{D}_\ell)^{I_K} \xrightarrow{\phi^{-1}_{K,\ell}} \Phi_{K,\ell},$$

where if $x \in A(K)_{\mathrm{tors},\ell}$, pick $\{x_i\}_{i=1}^\infty \in T_\ell$ such that $x = x_j$ for some $j \in \mathbb{N}$. Then set $g(x) :=$ class of $(\{x_i\}_{i=1}^\infty \otimes \ell^{-j})$. That the map $g$ is well defined and an isomorphism is proved in [Lor3], 3.4.

Let $\tau \in \Phi_{K,\ell}$. Let $y := \{y_i\}_{i=1}^\infty \otimes \ell^{-r} \in W_{L,\ell} \otimes \mathbb{Q}_\ell$ be such that $f_{W_{L,\ell}}(y) = \phi_{K,\ell}^{-1}(\tau)$. Let $\{t_i\}_{i=1}^\infty \otimes \ell^{-j}$ be such that $t = t_j$. Then in $(T_\ell \otimes \mathbb{D}_\ell)^{I_K}$, we have

$$\text{class of } (\{t_i\}_{i=1}^\infty \otimes \ell^{-j}) = \text{class of } (\{y_i\}_{i=1}^\infty \otimes \ell^{-r}).$$

Without loss of generality, we can assume that $j = r$. It is easy to check that the above equality in $T_\ell \otimes \mathbb{D}_\ell$ implies that $t_i = y_i$ for all $i = 1, \ldots, j$. Thus $t = t_j$ reduces in $\mathcal{T}_L$ since $W_{L,\ell}$ is canonically isomorphic to the Tate module of $\mathcal{T}_L$.

**6.8** Our next theorem describes an element $\tau \in \Phi(G)$ whose $\ell$-part belongs to the subgroup $\Theta_{K,\ell}^{[3]}$. To describe this element, we need to introduce the following notation. Let $(G, M, R)$ be any arithmetical graph. Let $(D, r)$ be a node of $G$. Let $(D_i, r_i)$, $i = 1, \ldots, d$, denote the vertices of $G$ linked to $D$. Assume that $(D_i \cdot D) = 1$ for all $i = 1, \ldots, d$, and that the numbering of the vertices $D_i$ is such that for $i = 1, \ldots, s$, the vertex $D_i$ belongs to a terminal chain $T_i$ attached at $D$, and for $i = s+1, \ldots, d$, the vertex $D_i$ is not on a terminal chain at $D$. We assume that $s \geq 2$. For simplicity, we will assume that $\gcd(r, r_i) = 1$, for all $i = 1, \ldots, s$. Thus the terminal vertex $C_i$ on $T_i$ has multiplicity 1. Let $\tau_i$ denote the image of $E(C_i, C_s)$ in $\Phi(G)$, $i = 1, \ldots, s-1$. Let

$$\tau := \sum_{i=1}^{s-1} r_i \tau_i.$$

To motivate this definition of $\tau$, let us note the following.

**Lemma 6.9** *Let $\ell$ be any prime. If $\mathrm{ord}_\ell(r) \leq \mathrm{ord}_\ell\left(\sum_{i=1}^s r_i\right)$, then $\langle \tau; \tau_i \rangle = 0$ in $\mathbb{Q}_\ell/\mathbb{Z}_\ell$, for all $i = 1, \ldots, s-1$.*

*Proof:* If $C$ is any vertex of $G$, let $r(C)$ denote its multiplicity. Then Lemma 2.8 shows the existence of integers $b_i$, $i = 1, \ldots, s$, such that

$$\sum_{\substack{C, C' \in T_i \\ (C \cdot C') = 1}} \frac{1}{r(C) r(C')} = \frac{b_i}{r \gcd(r, r_i)}.$$

Proposition 3.2 shows that

$$\langle \tau_i ; \tau_j \rangle = \begin{cases} \frac{b_s}{r} & \text{if } i \neq j \\ \frac{b_i}{r} + \frac{b_s}{r} & \text{if } i = j. \end{cases}$$

Thus, for $k = 1, \ldots, s-1$, we find that

$$\begin{aligned} \langle \tau ; \tau_k \rangle &= \sum_{i=1}^{s-1} r_i \langle \tau_i ; \tau_k \rangle = \left( \sum_{i=1}^{s-1} r_i \right) b_s / p_k + b_k r_k / r \\ &= \left( \sum_{i=1}^{s} r_i \right) b_s / r - b_s r_s / r + b_k r_k / r. \end{aligned}$$

Lemma 2.6 shows that $b_s r_s \equiv 1 \equiv b_k r_k \mod r$, and our hypothesis is that $\mathrm{ord}_\ell(r) \leq \mathrm{ord}_\ell \left( \sum_{i=1}^{s} r_i \right)$. Hence, $\langle \tau ; \tau_k \rangle = 0$ in $\mathbb{Q}_\ell / \mathbb{Z}_\ell$.

Assume that $(G, M, R)$ is associated to a curve $X/K$ whose jacobian has purely additive reduction. We have established in [Lor3], 3.13, that $\Theta^{[3]}_{K,\ell} = \Psi_{K,L,\ell} \cap \Psi^{\perp}_{K,L,\ell}$, where the orthogonal subgroup is computed with respect to the pairing 3.12 in [Lor3]. While no relationship between the pairing 3.12 and the pairing $\langle \ ; \ \rangle$ described in section 3 is fully established as of yet, one may certainly anticipate a relationship and, thus, we may expect that an element $\tau$ of $\Psi_{K,L,\ell}$ that is orthogonal to $\Psi_{K,L,\ell}$ under the pairing, $\langle \ , \ \rangle$ belongs to $\Theta^{[3]}_{K,\ell}$. Theorem 6.5 shows that the $\ell$-parts of $\tau_i$, $i = 1, \ldots, s-1$ and, thus, the $\ell$-part of $\tau$, belong to $\Psi_{K,L}$. Lemma 6.9 and Theorem 6.5 show that the $\ell$-part of $\tau$ is orthogonal to any element of $\Psi_{K,L}$ image of $E(C, C')$ with $\ell \nmid rr'$. Thus the $\ell$-part of $\tau$ is a 'good candidate' to be an element of $\Theta^{[3]}_{K,\ell}$, and Theorem 6.10 below describes some instances where the $\ell$-part of $\tau$ belongs to $\Theta^{[3]}_{K,\ell}$.

Note that if $s > 2$, then Theorem 5.4 shows that $\tau_i \neq 0$ for all $i = 1, \ldots, s-1$. But it may happen that $\tau$ is trivial in $\Phi(G)$, in which case the $\ell$-part of $\tau$ certainly belongs to $\Theta^{[3]}_{K,\ell}$, as in the following example (with $D$ being the node of multiplicity 4, and $\ell = 2$).

On the other hand, if $G$ contains a vertex $C$ with $\gcd(r, r(C)) = 1$ and $C \notin T_i$, for all $i = 1, \ldots, s$, then $\tau$ has order $r$ in $\Phi(G)$. Indeed, each $\tau_i$ has order $r$, thus the order of $\tau$ divides $r$. Let $\tau_C$ denote the image of $E(C_1, C)$ in $\Phi(G)$. Then $\langle \tau; \tau_C \rangle = b_1 r_C / r$. Thus, $r$ divides the order of $\tau$ (and of $\tau_C$).

**Theorem 6.10** *Let $X/K$ be a curve with a regular model $\mathcal{X}/\mathcal{O}_K$ and associated arithmetical graph $(G, M, R)$. Assume that the jacobian $A/K$ of $X/K$ has purely additive reduction over $\mathcal{O}_K$ and that the graph $G$ contains a node $(D, r)$ as in 6.8. Assume that $D$ and all components on the terminal chains attached to $D$ are smooth curves. Let $\ell \neq p$ be prime. Suppose that $r = \ell^{\mathrm{ord}_\ell(r)}$ and that $\mathrm{ord}_\ell(r_i) \geq \mathrm{ord}_\ell(r)$ for all $i = s+1, \ldots, d$. Then $\tau$ belongs to $\Theta^{[3]}_{K,\ell}$.*

*Proof:* Let $t_i \in A(K)_{\mathrm{tors},\ell}$ denote the unique torsion point in $A(K)$ whose image in $\Phi_{K,\ell}$ is equal to $\tau_i$. Let $P_i \in X(K)$ be such that $C_{P_i} = C_i$. Then $\pi_K(P_i - P_s) = \tau_i$. By hypothesis, the special fiber $\mathcal{A}_{K,k}$ is an extension of $\Phi_K$ by a unipotent group. Let $U_i$ denote the connected component of $\mathcal{A}_{K,k}$ such that $\pi_K(P_i - P_s) \in U_i$. Consider the natural map

$$\gamma_{K,L} : \mathcal{A}_K \times_{\mathrm{Spec}\mathcal{O}_K} \mathrm{Spec}\mathcal{O}_L \to \mathcal{A}_L.$$

Since $\mathcal{A}^0_{L,k}$ does not contain any unipotent group, the image of $U_i$ under $\gamma_{K,L}$ is a torsion point of $\mathcal{A}^0_{L,k}$ of order $r$. It follows that $\pi_L(P_i - P_s) = \pi_L(t_i)$. Hence, it follows from Lemma 6.7 that to prove Theorem 6.10, it is sufficient to exhibit an extension $F/K$ such that $\pi_F(\sum_{i=1}^{s-1} r_i(P_i - P_s)) \in \mathcal{T}_F$.

Let $M := M_r$. Consider the model $\mathcal{Y}/\mathcal{O}_M$ of $X_M/M$ associated as in 6.1 to $M/K$ and the model $\mathcal{X}/\mathcal{O}_K$. Let $E/k \subset \mathcal{Y}_k$ denote the strict transform of $D$ in $\mathcal{Y}$. It follows from 6.2 that if $P \in E$, then $P$ is a regular point on $\mathcal{Y}$. In particular, $E$ is a smooth curve. Let $b_{|E} : E \to D$ be the map obtaned by restriction from $\mathcal{Y} \to \mathcal{X}$. Let $k(D)$ be the function field of $D$. Choose a coordinate function $x$ in $k(D)$ so that when $D$ is identified with $\mathbb{A}^1/k \sqcup \{\infty\}$, then $\infty \neq D \cap D_i$, for all $i = 1, \ldots, s$. Let $a_i \in \mathbb{A}^1(k)$ denote the point $D \cap D_i$, $i = 1, \ldots, s$. The map $b_{|E} : E \to D$ is a cyclic Galois cover of degree $r$ ramified only above the points $a_i$, $i = 1, \ldots, s$ (we use here that $\mathrm{ord}_\ell(r_i) \geq \mathrm{ord}_\ell(r)$ for all $i = s+1, \ldots, d$). Thus

$$k(E) \cong k(x)[y]/(y^r - \prod_{i=1}^{s}(x - a_i)^{q_i})$$

for some positive integer $q_1, \ldots, q_s$ such that $r \mid \sum_{i=1}^{s} q_i$.

**Proposition 6.11** *We may choose $q_i = r_i$, for all $i = 1, \ldots, s$.*

Let $\xi$ denote a primitive $r$-th root of unity. Let $e_i$ denote the point of $E$ totally ramified above the point $a_i$ of $D$. Write $1 = \alpha_i r + \beta_i q_i$. Then, given any local uniformizer $\nu_i$ at $e_i$ with the property that $\sigma(\nu_i) = \xi^c \nu_i$ for some integer $c$, (such as $\nu_i = y^{\beta_i}(x-a_i)^{\alpha_i}$), we find that $c = \beta_i$ modulo $r$. In other words, $c$ is the inverse of $r_i$ modulo $r$.

Let us consider the map $c : \mathcal{X} \to \mathcal{X}'$, which contracts all components of $\mathcal{X}_k$ that belongs to a terminal chain attached to $D$. Thus $\mathcal{X}'$ is a normal model of $\mathcal{X}$ having exactly $s$ singular points $Q_1, \ldots, Q_s$ on $D$. Let $M := M_r$ and consider the base change maps $\mathcal{Y} \to \mathcal{X}$ and $\mathcal{Y}' \to \mathcal{X}'$, and the minimal desingularization maps $\mathcal{Z} \to \mathcal{Y}$ and $\mathcal{Z}' \to \mathcal{Y}'$. The map $c$ induces a map $c_{\mathcal{Y}} : \mathcal{Y} \to \mathcal{Y}'$. By construction, the multiplicity of $E$ in $\mathcal{Y}$ and $\mathcal{Z}$ is equal to 1. Thus every component $C$ of $\mathcal{Z}$ whose image in $\mathcal{X}$ belongs to a terminal chain attached to $D$ can be contracted by a map $\mathcal{Z} \to \mathcal{Z}''$ in such a way that the image of $C$ in $\mathcal{Z}''$ is a regular point of $\mathcal{Z}''$ (we use here the fact that all components of the terminal chains attached at $D$ are smooth). By minimality of the resolution, we find that we have maps $\mathcal{Z} \to \mathcal{Z}'' \to \mathcal{Z}'$. Thus every point of $\mathcal{Y}'$ in the image of $E$ is a regular point of $\mathcal{Y}'$.

Let us consider the action of the group $\mathrm{Gal}(M/K)$ on the scheme $\mathcal{Y}'$. The quotient of this action is the scheme $\mathcal{X}'$. Let $R_i$ denote the preimage of $Q_i$ in $\mathcal{Y}'$. As we mentioned above, $R_i$ is a regular point on $\mathcal{Y}'$. We may thus use the results on quotient singularities to describe the resolution of singularities at $Q_i$. Namely, the completed local ring at $R_i$ is of the form $\mathcal{O}_M[[z]]$, and $z$ can be chosen such that the action of $G$ on $\mathcal{O}_M[[z]]$ is linear: if $\sigma$ is a generator of $G$, then there exists a root of unity $\xi$ such that $\sigma(t_M) = \xi t_M$ and $\sigma(z) = \xi^{b_i} z$ for some $b_i \in \mathbb{N}$. (We use here the fact that the extension $M_r/K$ is tame.) Then the resolution of singularities at $Q_i$ is completely determined by the integer $b_i$. It follows (see for instance [Vie], 6.6) that in order to have a resolution of the type $\mathcal{X} \to \mathcal{X}'$, the integer $b_i$ must be congruent to the inverse of $r_i$ modulo $r$. Hence, we find that $q_i$ is congruent to $r_i$ modulo $r$.

**Lemma 6.12** *Let $s \geq 2$. Let $E/K$ denote the nonsingular complete model of the plane curve given by the equation*

$$y^r - \prod_{i=1}^{s}(x-a_i)^{r_i} = 0,$$

*with $a_i \in k$, $i = 1, \ldots, s$ and $\prod_{i \neq j}(a_i - a_j) \neq 0$. Assume that $\gcd(r, r_i) = 1$ for all $i = 1, \ldots, s$, and $r \mid \sum_{i=1}^{s} r_i$. Let $e_i$ denote the point of $E$ corresponding to the point $(a_i, 0)$. Then the divisor $e_i - e_s$ has order dividing $r$ in $\mathrm{Jac}(E)$, and $\sum_{i=1}^{s-1} r_i(e_i - e_s) = 0$ in $\mathrm{Jac}(E)$.*

*Proof:* The function $(x - a_i)/(x - a_s)$ belongs to the function field of $E$, and

$$\mathrm{div}((x - a_i)/(x - a_s)) = r(e_i - e_s).$$

Moreover, let $d := \sum_{i=1}^{s} r_i/r$. Then $(y/(x - a_s)^d)^r = \prod_{i=1}^{s-1} \left(\frac{x-a_i}{x-a_s}\right)^{r_i}$. Thus, $\sum_{i=1}^{s-1} r_i(e_i - e_s) = \mathrm{div}(y/(x - a_s)^d)$.

Let us now conclude the proof of Theorem 6.10. We are going to show that $\pi_{ML}(\sum_{i=1}^{s-1} r_i(P_i - P_s)) \in \mathcal{T}_{ML}$. The special fiber of the Néron model $\mathcal{A}_M^0$ is isomorphic to $\mathrm{Pic}^0(\mathcal{Z}_k/k)$. The group scheme $\mathrm{Pic}^0(\mathcal{Z}_k/k)$ is an extension of the abelian variety $B_M := \prod_{C \subseteq \mathcal{Z}_k} \mathrm{Jac}(C)$ by the product of a unipotent group $U_M$ and a torus $\mathcal{T}_M$. Lemma 6.12 implies that the image of $\pi_M(\sum_{i=1}^{s-1} r_i(P_i - P_s))$ in $B_M$ is trivial. Thus, the image of $\pi_{ML}(\sum_{i=1}^{s-1} r_i(P_i - P_s))$ in $B_{ML}$ is also trivial. Since $A_{ML}/ML$ has semi-stable reduction, we find that $\pi_{ML}(\sum_{i=1}^{s-1} r_i(P_i - P_s))$ belongs to $\mathcal{T}_{ML}$.

# 7 Partial converses for Theorem 6.5

Let $X/K$ be a curve. Let $\mathcal{X}/\mathcal{O}_K$ be a regular model of $X/K$. Let $(C, r)$ and $(C', r')$ be two distinct components of $\mathcal{X}_k$. Let $\ell \neq p$ be a prime. In view of Theorem 6.5, it is natural to wonder whether it is true that if the $\ell$-part of $E(C, C')$ belongs to $\Psi_{K,L}$, then the pair $(C, C')$ is weakly connected and $\ell$-breakable. Let us make the following conjectures.

**Conjecture 7.1** *Let $\ell \neq p$. Let $(C, C')$ be a weakly connected pair such that $\ell \nmid rr'$. If $(C, C')$ is not $\ell$-breakable, then the $\ell$-part of $E(C, C')$ does not belong to $\Psi_{K,L}$.*

If $\ell = p$, Conjecture 7.1 does not hold, as can be seen on the following example with $p = 2$. Consider an elliptic curve $X/K$ with reduction $I_\nu^*$, $\nu > 1$, and with potentially good reduction. Then the graph of the reduction of $X$ contains pairs that are not $p$-breakable. On the other hand, since $X$ has potentially good reduction, $\Psi_{K,L} = \Phi_K$.

Clearly, Conjecture 7.1 implies that the $\ell$-part of $E(C, C')$ is not trivial in $\Phi_K$. This fact is proved in 3.3. A second partial converse to Theorem 6.5 could be the following:

**Conjecture 7.2** Let $\ell \neq p$. Let $(C, C')$ be a multiply connected pair such that $r = r' = 1$. If the $\ell$-part of $E(C, C')$ is not trivial in $\Phi_K$, then the $\ell$-part of $E(C, C')$ does not belong to $\Psi_{K,L}$.

It is possible that Conjecture 7.2 holds more generally for $\ell$ any prime, and $\gcd(\ell, rr') = 1$. As evidence for these two conjectures, we offer the following two theorems.

**Theorem 7.3** *Let $\ell \neq p$. Let $\mathcal{X}/\mathcal{O}_K$ be a regular model of a curve $X/K$, with associated arithmetical graph $(G, M, R)$. Let $(C, C')$ be a weakly connected pair with $r = r' = 1$. If $(C, C')$ is not $\ell$-breakable in $G$, then the image $\tau$ of $E(C, C')$ does not belong to $\Psi_{K,L}$.*

*Proof:* Lemma 6.4 shows that we only need to exhibit an extension $M/K$ such that $\tau \in \Psi_{K,M}$. Consider the base change $M_\ell/K$ and the associated map $\rho \circ b : \mathcal{Z} \to \mathcal{X}$ introduced in 6.1. We denote again by $C$ and $C'$ the preimages in $\mathcal{Z}$ of $C$ and $C'$. Lemma 6.3 shows that the pair $(C, C')$ is multiply connected in the graph of the model $\mathcal{Z}$. Thus Theorem 7.3 follows from our next theorem.

**Theorem 7.4** *Let $\mathcal{X}/\mathcal{O}_K$ be a regular model of a curve $X/K$, with associated arithmetical graph $(G, M, R)$. Let $(C, C')$ be a multiply connected pair with $r = r' = 1$. Then $E(C, C')$ does not belong to $\Psi_{K,L}$.*

*Proof:* Consider the base change $L/K$ and the associated map $b \circ \rho : \mathcal{Z} \to \mathcal{X}$. We denote again by $C$ and $C'$ the preimages in $\mathcal{Z}$ of $C$ and $C'$. It is easy to check (6.3) that the pair $(C, C')$ is multiply connected in $\mathcal{Z}$. Let $\mathcal{Z}_{min}$ denote the minimal regular model of $X_L/L$. Since $(\mathcal{Z}_{min})_k$ is a reduced curve, its associated graph is reduced (see 2.1). The contraction map $\gamma : \mathcal{Z} \to \mathcal{Z}_{min}$ may contract connected curves. However, the dual graph of each configuration that may be contracted is a tree. Moreover, because $C$ and $C'$ intersect at least two components of $\mathcal{Z}$ and have both multiplicity one, neither $C$ nor $C'$ can be contracted by $\gamma$. Thus the image of $C$ and $C'$ in $\mathcal{Z}_{min}$ is a multiply connected pair. To conclude the proof of Theorem 7.4, we use Corollary 2.3 in [Lor4] (see 2.1), which states that in a reduced graph $G$, a pair of vertices is trivial in $\Phi(G)$ if and only if the pair is weakly connected.

**Corollary 7.5** *Let $(G, M, R)$ be any arithmetical graph. If $G$ contains a multiply connected pair $(C, C')$, with $r = r' = 1$, then $|\Phi(G)| \neq 1$.*

*Proof:* Winters' Existence Theorem [Win] implies the existence of a field $F$ with a discrete valuation of equicharateristic 0, and a smooth and proper curve $Y/F$ having a model over $\mathcal{O}_F$ whose associated arithmetical graph is the given graph $(G, M, R)$. Apply 7.4.

It would be interesting to find a direct proof of 7.5 that does not rely on the theory of degenerations of curves. Example 7.9 below shows that if only one vertex of $G$ has multiplicity equal to 1, then it may happen that $|\Phi(G)| = 1$. Let us now show with the help of some examples that the above conjectures do not hold in general if $\ell$ divides $rr'$.

**Example 7.6** Let $a$ and $b$ be two positive integers. Consider the arithmetical graph $G$ given by:

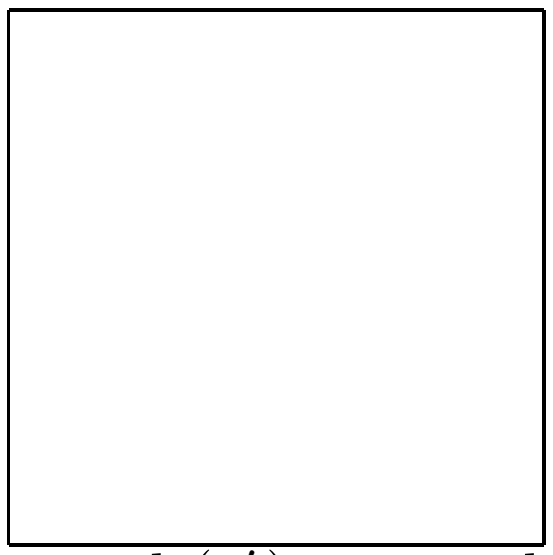

where $\mathrm{ord}_\ell(r) = a$, $\mathrm{ord}_\ell(r_1) = \mathrm{ord}_\ell(r'_1) = 1$, $\mathrm{ord}_\ell(s) = a + b$, and $\mathrm{ord}_\ell(s_1) = \mathrm{ord}_\ell(s'_1) = 1$. Note that the relations $t_1$ divides $r + t_2$, $t_2$ divides $t_1 + t_3$, ..., $t_k$ divides $t_{k-1} + s$, show that the sequence

$$s\ell^{-a}, t_k\ell^{-a}, \ldots, t_1\ell^{-a}, r\ell^{-a},$$

can be continued using Euclid's algorithm with $t_1\ell^{-a}$ and $r\ell^{-a}$ into a sequence

$$S := \{s\ell^{-a}, t_k\ell^{-a}, \ldots, t_1\ell^{-a}, r\ell^{-a}, u_1, \ldots, u_w\},$$

as in [Lor2], 2.4, and that this sequence $S$ can be considered as the sequence of multiplicities of a terminal chain of some arithmetical graph. In particular, Lemma 2.8 shows that

$$\frac{1}{u_w u_{w-1}} + \ldots + \frac{1}{u_2 u_1} + \frac{1}{u_1 r\ell^{-a}} + \frac{1}{r\ell^{-a} t_1 \ell^{-a}} + \cdots + \frac{1}{t_k \ell^{-a} s\ell^{-a}} = \frac{\beta}{s\ell^{-a} u_w}$$

for some integer $\beta$ coprime to $s\ell^{-a}$, and

$$\frac{1}{u_w u_{w-1}} + \ldots + \frac{1}{u_2 u_1} + \frac{1}{u_1 r\ell^{-a}} = \frac{\gamma}{r\ell^{-a}u_w}$$

for some integer $\gamma$ coprime to $r\ell^{-a}$,

**Lemma 7.7** *The $\ell$-part of the group $\Phi(G)$ is cyclic of order $\ell^{2a+b}$ and is generated by the $\ell$-part of the image of $E(B,C)$.*

*Proof:* Proposition 9.6/6 of [BLR] shows that $|\Phi(G)| = rs/r_n r'_{n'} s_m s'_{m'}$, so its $\ell$-part has order $\ell^{2a+b}$. Let $\tau \in \Phi(G)$ denote the image of $E(B,C)$. Consider the pairing $\langle\ ,\ \rangle$ introduced in 3.1. To show that the $\ell$-part of $\tau$ has order $\ell^{2a+b}$, it is sufficient to show that the order of $\langle\tau,\tau\rangle$ in $\mathbb{Q}/\mathbb{Z}$ is divisible by $\ell^{2a+b}$. We compute this order using Proposition 3.2:

$$\langle\tau,\tau\rangle = \left(\frac{1}{r_n r_{n-1}} + \cdots + \frac{1}{r_1 r}\right) + \left(\frac{1}{rt_1} + \cdots + \frac{1}{t_{k-1}t_k} + \frac{1}{t_k s}\right) + \left(\frac{1}{ss_1} + \cdots + \frac{1}{s_{m-1}s_m}\right).$$

Lemma 2.8 implies that there exist integers $c_n$ and $d_m$ coprime to $\ell$ such that

$$\langle\tau,\tau\rangle = \frac{c_n}{r} + \left(\frac{\beta\ell^{-a}}{su_w} - \frac{\gamma\ell^{-a}}{ru_w}\right) + \frac{d_m}{s}.$$

Regarded as elements in $\mathbb{Q}_\ell/\mathbb{Z}_\ell$, $c_n/r$ has order $\ell^a$, $d_m/s$ has order $\ell^{a+b}$, $\gamma/\ell^a r u_w$ has order at most $\ell^{2a}$. Since $b > 0$, $\ell \nmid \beta$ and, thus, $\beta/\ell^a s u_w$ has order $\ell^{2a+b}$. Since $a > 0$, we find that $\langle\tau,\tau\rangle$ has order $\ell^{2a+b}$ in $\mathbb{Q}_\ell/\mathbb{Z}_\ell$, which concludes the proof of Lemma 7.7.

Let us consider now the case where $r = \ell^a$ and $s = \ell^{a+b}$. Proposition 4.3 shows that the order of $(A,B)$ is $\ell^a$, while the order of $(C,D)$ is $\ell^{a+b}$. Let $C_r$ and $C_s$ denote the nodes of multiplicity $r$ and $s$, respectively. The pair $(C_r, C_s)$ is weakly connected but not $\ell$-breakable. We find using Remark 3.5 that $\ell^{a+b}E(B,C) = E(C_r, C_s)$. Winters' Existence Theorem [Win] implies the existence of a field, say $K$, with a discrete valuation of equicharateristic 0, and a smooth and proper curve $X/K$ having a model over $\mathcal{O}_K$ whose associated arithmetical graph is the given graph $(G, M, R)$. Since $\ell$ is not equal to the residue characteristic of $K$, Theorem 6.5 shows that $E(C,D)$ belongs to $\Psi_{K,L}(\mathrm{Jac}(X))$. Thus $\ell^{a+b}$ divides $|\Psi_{K,L}|$. Hence, the group $\Phi_K/\Psi_{K,L}$ is killed by $\ell^a$ and any element of $\Phi_K$ which is $\ell^a$-divisible in $\Phi_K$ belongs to $\Psi_{K,L}$. It follows that the pair $(C_r, C_s)$ belongs to $\Psi_{K,L}$ but is not $\ell$-breakable.

Theorem 6.5 shows that $E(A,B)$ belongs to $\Psi_{K,L}$. Thus, since $\Phi(G)$ is cyclic, the element $E(A,B)$ is a multiple of $E(C,D)$, and is divisible by $\ell$ in $\Psi_{K,L}$. We shall see in the next section that this phenomenon cannot occur if $\mathrm{Jac}(X)$ has potentially good $\ell$-reduction.

Let $L/K$ denote the extension minimal with the property that $X_L/L$ has semistable reduction. Let $t_L$ and $a_L$ denote the toric and abelian ranks of $\mathrm{Jac}(X_L)/L$, respectively. When $\ell$ is not the residue characteristic, one can show that $t_L = \ell^a - 1$, and $a_L = (\ell^{a+b} - \ell^a)/2$. It is shown in [Lor3], 1.7, (using the fact that $\Phi$ is cyclic) that $|\Psi_{K,L}| - 1 \leq 2a_L + t_L$. It follows from this bound that $|\Psi_{K,L}| = \ell^{a+b}$, and $|\Psi_{K,L}| - 1 = 2a_L + t_L$. It would be very interesting to know what are the possible values of the integers $t_L$ and $a_L$ when $\ell$ is the residue characteristic of a field $K$ and there exists a curve $X/K$ having a model over $\mathcal{O}_K$ whose associated arithmetical graph is the given graph $(G, M, R)$. We conjecture that in this case $t_L \leq \ell^a - 1$.

Note that when $\Phi_{K,\ell}$ is cyclic, the subgroup $\Theta^{[3]}_{K,\ell}$ is completely determined by $\Psi_{K,L,\ell}$. Indeed, $\Phi_{K,\ell}$ is endowed with a perfect pairing such that $\Theta^{[3]}_{K,\ell}$ is the orthogonal of $\Psi_{K,L,\ell}$. Thus, if $|\Phi_{K,\ell}| = \ell^{2a+b}$ and $|\Psi_{K,L,\ell}| = \ell^{a+b}$, then $|\Theta^{[3]}_{K,\ell}| = \ell^a$.

**Remark 7.8** Let us use a graph $G$ of the type introduced in Example 7.6 to exhibit an example where the $\ell$-part of the group $\Phi(G)$ is not generated by the images of the elements $E(C,C')$, with $\gcd(\ell, rr') = 1$. The multiplicities of $G$ are as follows. Let $r_n = r_1 = 1$, $r = 4$, $r'_{n'} = r'_1 = 1$, $t_1 = 8$, $t_2 = 12$, $s = 16$, $s_1 = s'_1 = 10$, $s_2 = s'_2 = 4$, and $s_3 = s'_3 = 2$. The order of $\Phi(G)$ is equal to 16. Let then $\ell = 2$. The pair $(A,B)$ has order 2 and is the only pair with $2 \nmid rr'$. The image $\tau$ of $E(B,C)$ has order 16 since $\langle \tau, \tau \rangle$ is easily computed to have order at least 16.

**Example 7.9** Let $(C,C')$ be a multiply connected pair. It is not true in general that the $\ell$-part of $E(C,C')$ is not trivial in $\Phi_K/\Psi_{K,L}$. Indeed, it may happen that $E(C,C')$ belongs to $\Psi_{K,L}$ simply because $\Phi_K$ itself is trivial. Consider the following example:

The order of $\Phi_K$ equals $2/\gcd(r, r-2)$. Hence, when $r$ is even, $|\Phi_K| = 1$. When $|\Phi_K| = 2$, $E(C, D)$ is a generator.

# 8 The case of potentially good $\ell$-reduction

Our goal in this section is to prove the following theorem. Recall the definitions introduced in 1.4.

**Lemma 8.1** *Let $A/K$ be a principally polarized abelian variety. Let $\ell$ be a prime, $\ell \neq p$. Assume that $A/K$ has potentially $\ell$-good reduction. Then $\Psi_{K,L,\ell} = (\Phi_K)_\ell$.*

*Proof:* It is shown in [Lor3] (see 3.22, with 3.21 (ii) and 2.15 (ii)), that the kernel of the map $\Phi_{K_\ell} \to \Phi_L$ is killed by $[L : K_\ell]$. Thus $\Psi_{K_\ell,L,\ell} = (0)$. It also follows from [Lor3], using the fact that $t_{K_\ell} = 0$, that $\Psi_{K_\ell,L,\ell} = (\Phi_{K_\ell})_\ell$. Thus, since $(\Phi_{K_\ell})_\ell = (0)$, we find that $\Psi_{K,L,\ell} = (\Phi_K)_\ell$.

**Theorem 8.2** *Let $X/K$ be a curve. Let $\ell$ be a prime, $\ell \neq p$, and assume that* $\mathrm{Jac}(X)/K$ *has potentially good $\ell$-reduction. Let $P, Q \in X(K)$ with $C_P \neq C_Q$. Then the $\ell$-part $\tau_\ell$ of the image of $P - Q$ in $\Phi_K$ belongs to $\Psi_{K,L,\ell}$. If $\tau_\ell$ is not trivial, then it is not $\ell$-divisible in $\Phi_K$.*

*Proof:* Lemma 8.1 shows that $\tau_\ell$ belongs to $\Psi_{K,L}$. Theorem 8.2 is a consequence of Theorem 8.3 below, which pertains only to arithmetical graphs. Indeed, Proposition 1.7 in [Lor2] shows that if $\mathrm{Jac}(X)/K$ has potentially good $\ell$-reduction, then there exists a model $\mathcal{X}/\mathcal{O}_K$ of $X/K$ whose associated graph $G$ is a tree satisfying Condition $C_\ell$ stated in 1.5 of [Lor2].

**Theorem 8.3** *Let $(G, M, R)$ be an arithmetical tree. Let $\ell$ be any prime. Let $(C, r)$ and $(C', r')$ be two vertices of $G$ such that $\ell \nmid rr'$. If $G$ satisfies Condition $C_\ell$, then the $\ell$-part of $E(C, C')$ has order $\lambda(C, C')$. Moreover, if the $\ell$-part of $E(C, C')$ is not trivial, then it is not $\ell$-divisible.*

*Proof:* Since $G$ is a tree, every pair $(C, C')$ is weakly connected. Condition $C_\ell$ implies that any two vertices $C$ and $C'$ with $\ell \nmid rr'$ form an $\ell$-breakable weakly connected pair. Thus we may use Theorem 5.4 to compute the order of $E(C, C')$. Let us now show that $E(C, C')$ is not divisible by $\ell$ if it is not trivial. If the path $\mathcal{P}$ connecting $C$ to $C'$ does not contain any node, then Theorem 5.4 shows that the $\ell$-part of the order of $E(C, C')$ is trivial and, thus, in this case the statement of Theorem 8.3 does not apply. Let us now assume that $\mathcal{P}$ contains at least one node.

We claim that Theorem 8.3 holds if it holds in the special case where $\mathcal{P}$ has only one node. Indeed, if the path $\mathcal{P}$ connecting $C$ to $C'$ contains more than one node, use Proposition 5.3 to break the tree $G$ into several trees $G_1, \ldots, G_m$, each having a weakly connected $\ell$-breakable pair of terminal vertices $C_i$ and $C'_i$ connected by a path having at most one node. Each tree $G_j$ satisfies Condition $C_\ell$. The construction of the graphs $G_i$ is such that $\Phi_\ell(G) \cong \prod_{i=1}^m \Phi_\ell(G_j)$, Moreover, the image of the $\ell$-part of $E(C, C')$ in $\Phi(G_j)$ is the $\ell$-part of $E(C_i, C'_i)$. Thus, the $\ell$-part of $E(C, C')$ is not $\ell$-divisible in $\Phi(G)$ if and only if the $\ell$-part of $E(C_i, C'_i)$ is not $\ell$-divisible in $\Phi(G_i)$ for some $i$.

Consider now the case where $(G, M, R)$ is an arithmetical tree satisfying Condition $C_\ell$, with a pair of terminal vertices $(C, r)$ and $(C', r')$ such that $\ell \nmid rr'$, and such that the path $\mathcal{P}$ connecting $C$ to $C'$ in $G$ contains a unique node $(D, r_D)$. Let $\nu$ denote the total number of nodes of $G$. We proceed by induction on $\nu$. Assume that $(G, M, R)$ is an arithmetical tree with only one node $(D, r)$. Let $(C_1, r_1), \ldots, (C_d, r_d)$ denote the vertices of $G$ adjacent to $D$. The vertices $(D, r)$ and $(C_i, r_i)$ are on a unique terminal chain $T_i$, with terminal vertex of multiplicity $s_i := \gcd(r, r_i)$. We may always order the vertices $C_i$ such that

$$\mathrm{ord}_\ell(s_1) \geq \cdots \geq \mathrm{ord}_\ell(s_{d-1}) = \mathrm{ord}_\ell(s_d) = 1$$

(see [Lor2], 2.7). In particular, $\ell \nmid s_{d-1}s_d$. Denote by $(D_i, s_i)$ the terminal vertex of the chain $T_i$. Without loss of generality, we may assume that $C = D_d$ and $C' = D_{d-1}$. It is shown in [Lor2], 2.1, that the group $\Phi_\ell(G)$ is isomorphic to $\prod_{i=1}^{d-2} \mathbb{Z}/\ell^{\mathrm{ord}_\ell(r/s_i)}\mathbb{Z}$. It follows from Proposition 4.3 that the $\ell$-part of $E(D_{d-1}, D_d)$ has order $\ell^{\mathrm{ord}_\ell(r/s_{d-2})}$ in $\Phi(G)$. Thus, the $\ell$-part of $E(D_{d-1}, D_d)$ is not divisible by $\ell$ in $\Phi(G)$.

Consider now the case where $\nu > 1$ and proceed as follows. Pick an edge $e$ of $G$ such that one of the two components of $G \setminus \{e\}$ contains a single node $B$, with $B \neq D$. The component that does not contain $B$ can be completed into a new arithmetical graph $G'$, as in [Lor2], page 165. The graph $G'$ satisfies Condition $C_\ell$, and has $\nu - 1$ nodes. Thus we may apply the induction hypothesis and obtain that $E(C, C')$ is not divisible by $\ell$ in $\Phi(G')$. The discussion on page 165 of [Lor2] shows that $\Phi_\ell(G)$ contains $\Phi_\ell(G')$ as a direct summand. Since $E(C, C')$ is not divisible by $\ell$ in $\Phi(G')$, $E(C, C')$ is not divisible by $\ell$ in $\Phi(G)$.

**Remark 8.4** The fact in Theorem 8.2 that a point of the form $P - Q$ is not divisible by $\ell$ does not generalize to a statement pertaining to the group $\Psi_{K,L}$. Indeed, when $\Psi_{K,L} \neq \Phi_K$, Example 7.6 exhibits a point of the form $P - Q$ in $\Psi_{K,L}$ that is divisible by $\ell$ in $\Psi_{K,L}$.